\documentclass[12pt]{article}
\usepackage{amsmath,amsthm,amsfonts,amssymb}
\usepackage{fullpage,microtype}
\usepackage[authoryear]{natbib}

\usepackage{pgf}
\usepackage{subcaption}
\usepackage[margin=3em,font=small,labelfont=bf]{caption}

\usepackage{algorithm}
\usepackage{algorithmic}

\usepackage{xcolor}
\newcommand{\bo}[1]{}

\newtheorem{theorem}{Theorem}[section]

\newtheorem{claim}[theorem]{Claim}
\newtheorem{observation}[theorem]{Observation}
\newtheorem{lemma}[theorem]{Lemma}
\newtheorem{corollary}[theorem]{Corollary}
\newtheorem{proposition}[theorem]{Proposition}

\theoremstyle{definition} \newtheorem{definition}[theorem]{Definition}
\theoremstyle{definition}

\newenvironment{myproof}{\small \paragraph{\small Proof.}}{\hfill \qed \endproof \vskip1em}

\let\originalleft\left
\let\originalright\right
\renewcommand{\left}{\mathopen{}\mathclose\bgroup\originalleft}
\renewcommand{\right}{\aftergroup\egroup\originalright}

\newcommand{\reals}{\mathbb{R}}
\newcommand{\extreals}{\overline{\reals}}

\newcommand{\cl}[1]{\textrm{cl}\left(#1\right)} 
\newcommand{\inter}[1]{\textrm{int}\left(#1\right)} 
\newcommand{\effdom}[1]{\textrm{effdom}\left(#1\right)} 

\newcommand{\Supp}[1]{\textrm{Supp}\left(#1\right)} 

\newcommand{\sign}{\textrm{sign}}  
\newcommand{\st}{~:~}  

\renewcommand{\P}{\mathcal{P}}  
\newcommand{\Y}{\mathcal{Y}}  

\begin{document}
\title{Linear Functions to the Extended Reals\footnote{Many thanks to Rafael Frongillo for discussions and suggestions; and to Terry Rockafellar, an anonymous reviewer, and the authors of \citet{dudik2022convex} for comments.}}
\author{Bo Waggoner  \\ {\small University of Colorado, Boulder}}
\date{\today}
\maketitle

\begin{abstract}
  This paper investigates functions from $\reals^d$ to $\reals \cup \{\pm \infty\}$ that satisfy axioms of linearity wherever allowed by extended-value arithmetic.
  They have a nontrivial structure defined inductively on $d$, and unlike finite linear functions, they require $\Omega(d^2)$ parameters to uniquely identify.
  In particular they can capture vertical tangent planes to epigraphs: a function is convex if and only if it is the pointwise supremum of a family of ``affine extended'' functions.
  These results are applied to the well-known characterization of proper scoring rules for the finite-dimensional case: it is carefully and rigorously extended here to a more constructive form.
  In particular it is investigated when proper scoring rules can be constructed from a given convex function.
\end{abstract}
\vskip2em

\section{Introduction}
The extended real number line, denoted $\extreals = \reals \cup \{\pm \infty\}$, is widely useful particularly in convex analysis.
But I am not aware of an answer to the question: What would it mean to have a ``linear'' $f: \reals^d \to \extreals$?
The extended reals have enough structure to hope for a useful answer, but differ enough from a vector space to need investigation.
A natural approach is that $f$ must satisfy the usual linearity axioms of homogeneity and additivity whenever legal under extended-reals arithmetic.
Here is an example when $d=3$:
  \[ f(x,y,z) = \begin{cases}  \text{if $z > 0$} &\rightarrow~~ \infty  \\
                               \text{if $z < 0$} &\rightarrow~~ -\infty  \\
                               \text{if $z = 0$} &\rightarrow~~ 
        \begin{cases} \text{if $y > 0$} &\rightarrow~~ \infty  \\
                      \text{if $y < 0$} &\rightarrow~~ -\infty  \\
                      \text{if $y = 0$} &\rightarrow~~ x . \end{cases} \end{cases} \]

We will see that all \emph{linear extended} functions (Definition \ref{def:linext}) on $\reals^d$ have the above inductive format, descending dimension-by-dimension until reaching a finite linear function (Proposition \ref{prop:alg-correct}).
In fact, a natural representation of this procedure is parsimonious, requiring as many as $\Omega(d^2)$ real-valued parameters to uniquely identify an extended linear function on $\reals^d$ (Proposition \ref{prop:parsim}).
This structure raises the possibility, left for future work, that linear extended functions are significantly nontrivial on infinite-dimensional spaces.

Linear extended functions arise naturally as \emph{extended subgradients} (Definition \ref{def:extsub}) of a convex function $g$.
These can be used to construct \emph{affine extended} functions (Definition \ref{def:affext}) capturing vertical supporting hyperplanes to convex epigraphs.
Along such a hyperplane, intersected with the boundary of the domain, $g$ can sometimes have the structure of an arbitrary convex function in $d-1$ dimensions.
For example, the above $f$ is an extended subgradient, at the point $(1,0,0)$, of the discontinuous convex function
  \[ g(x,y,z) = \begin{cases}  \text{if $z > 0$} &\rightarrow~~ \infty  \\
                               \text{if $z < 0$} &\rightarrow~~ 0   \\
                               \text{if $z = 0$} &\rightarrow~~
        \begin{cases} \text{if $y > 0$} &\rightarrow~~ \infty  \\
                      \text{if $y < 0$} &\rightarrow~~ 0  \\
                      \text{if $y = 0$} &\rightarrow~~ \frac{1}{2}x^2 . \end{cases} \end{cases} \]
A function is convex if and only if it is the pointwise supremum of affine extended functions (Proposition \ref{prop:sup-char}).
We also find it is convex if and only if it has an extended subgradient at every point where it is finite, assuming a convex effective domain (Proposition \ref{prop:subgrad-char}).

\vskip1em
The role of vertical tangents to epigraphs has certainly been noted and appreciated in convex analysis.
For example, \citet{rockafellar1998variational} utilizes ``horizon subgradients''.
However, such approaches motivated by e.g. optimization problems focus on niceness properties, e.g. benign topological properties, not explored here.
Recently, \citet{dudik2022convex} independently discover linear extended functions with an optimization and topological motivation.
Beginning with an algorithmic definition, essentially Algorithm \ref{alg:linext}, \cite{dudik2022convex} shows that linear extended functions are characterized as being both convex and concave while vanishing at the origin.
That work also considers extended subgradients, there called \emph{astral dual subgradients}.

\vskip1em
\paragraph{Proper scoring rules.}
The original motivation for this investigation was the study of scoring rules: functions $S$ assigning a score $S(p,y)$ to any prediction $p$ (a probability distribution over outcomes) and observed outcome $y$.
This $S$ is called \emph{proper} if reporting the true distribution of the outcome maximizes expected score.
A well-known characterization states that proper scoring rules arise as and only as subgradients or ``subtangents'' of convex functions~\citep{mccarthy1956measures,savage1971elicitation,schervish1989general}.

Modern works~\citep{gneiting2007strictly,frongillo2021general} generalize this characterization to allow scores of $-\infty$, corresponding to convex functions that are not subdifferentiable\footnote{They also consider infinite-dimensional outcome spaces, which are not treated here.}.
These works therefore replace subgradients with briefly-defined generalized ``subtangents'', but they slightly sidestep the questions tackled head-on here: existence of these objects and their roles in convex analysis.
So these characterizations are not fully constructive.
In particular, it may surprise even experts that rigorous answers are not available to the following questions, given a convex function $g$: \emph{(a)} under what conditions can one construct a proper scoring rule from it? \emph{(b)} when is the resulting scoring rule necessarily \emph{strictly} proper? \emph{(c)} do the answers depend on which subgradients of $g$ are used, when multiple choices are available?

Section \ref{sec:proper-scoring-rules} uses the machinery of extended linear functions to prove a ``construction of proper scoring rules'' that answers such questions for the finite-outcome case.

First, Theorem \ref{thm:proper-char} allows predictions to range over the entire simplex.
It implies a slightly-informal claim of \citet{gneiting2007strictly}, Section 3.1: \emph{Any convex function on the probability simplex gives rise to a proper scoring rule (using any choices of its extended subgradients).}
Furthermore, it shows that a strictly convex function can \emph{only} give rise to strictly proper scoring rules and vice versa.
These facts are likely known by experts in the community; however, I do not know of formal claims and proofs.
This may be because, when scoring rules can be $-\infty$, proofs seem to require significant formalization and investigation of ``subtangents''.
In this paper, this investigation is supplied by linear extended functions and characterizations of (strictly) convex functions as those that have (uniquely supporting) extended subgradients everywhere.

Next, Theorem \ref{thm:proper-subset-char} allows prediction spaces $\P$ to be any nonempty subset of the simplex.
It sharpens the characterizations of \citet{gneiting2007strictly,frongillo2021general}, for the finite-outcome case, by showing $S$ to be proper if and only if it can be constructed from extended subgradients of a convex $g$ that is \emph{interior-locally-Lipschitz} (Definition \ref{def:inter-local}).
It also answers the construction questions \emph{(a)-(c)} above, e.g. showing that given such a $g$, some but not necessarily all selections of its extended subgradients give rise to proper scoring rules.
An important remaining direction is an analogous construction for truthful scoring functions on general type spaces, as studied by \citet{frongillo2021general}.

\paragraph{Preliminaries.}
Functions in this work are defined on Euclidean spaces of dimension $d\geq 0$.
I ask the reader's pardon for abusing notation slightly: I use $\reals^d$ refer to any such space.
For example, I may argue that a certain function exists on domain $\reals^d$, then refer to that function as being defined on a given $d$-dimensional subspace of $\reals^{d+1}$.

Let $g: \reals^d \to \extreals$.
The \emph{effective domain} $\effdom{g}$ of $g$ is $\{x \in \reals^d \st g(x) \neq \infty\}$.
The function $g$ is \emph{convex} if its \emph{epigraph} $\{(x,y) \in \reals^d \times \reals : y \geq g(x)\}$ is a convex set.
Equivalently, it is convex if for all $x,x' \in \effdom{g}$ and all $0 < \rho < 1$, $g(\rho \cdot x + (1-\rho)x') \leq \rho g(x) + (1-\rho) g(x')$, observing that this sum may contain $-\infty$ but not $+\infty$.
It is \emph{strictly} convex on a convex set $\P \subseteq \effdom{g}$ if the previous inequality is always strict for $x,x' \in \P$ with $x \neq x'$.
A function $g$ is \emph{concave} if $-g$, i.e. the function $x \mapsto -g(x)$, is convex.

We say a function $h$ \emph{minorizes} $g$ if $h(x) \leq g(x)$ for all $x$.
$\cl{A}$ denotes the closure of the set $A$ and $\inter{A}$ its interior.
For $a \in \extreals$, the sign function is $\sign(a) = 1$ if $a > 0$, $\sign(a) = 0$ if $a=0$, and $\sign(a) = -1$ if $a < 0$.


\paragraph{The extended reals.}
The extended reals, $\extreals = \reals \cup \{\pm \infty\}$, have the following rules of arithmetic.
For $a,b \in \extreals$, the sum $a + b$ is illegal if $\{a,b\} = \{\pm \infty\}$.
Otherwise it is legal.
We have $\infty + \infty = \infty$ and analogously for $-\infty$.
For any $\alpha \in \reals$, $\alpha + \infty = \infty$, $\alpha - \infty = -\infty$, and
  \[ \alpha \cdot \infty = \begin{cases} \infty & \alpha > 0  \\ 0  & \alpha = 0  \\ -\infty  & \alpha < 0 . \end{cases} \]
Addition is associative and commutative as long as it is legal; multiplication of multiple reals and possibly one non-real is associative and commutative.
Multiplication by a real distributes over legal sums.

$\extreals$ has the following rules of comparison: $-\infty < \beta < \infty$ for every $\beta \in \reals$.
The supremum of a subset of $\extreals$ is the least element of $\extreals$ that upper-bounds every member of the subset.
In particular, $\sup A = \infty$ if $\infty \in A$ or if $A$ contains an unbounded-above set of reals.
The analogous facts hold for the infimum.
Finally, $\inf \emptyset = \infty$ and $\sup \emptyset = -\infty$.
I will not put a topology on $\extreals$ in this work.

\section{Linear extended functions}
The following definition makes sense with a general real vector space $\mathcal{X}$ in place of $\reals^d$, but it remains to be seen which results can extend.

\begin{definition} \label{def:linext}
  Call the function $f: \reals^d \to \extreals$ a \emph{linear extended function} if:
  \begin{enumerate}
    \item (scaling) For all $x \in \reals^d$ and all $\alpha \in \reals$: $f(\alpha x) = \alpha f(x)$.
    \item (additivity) For all $x,x' \in \reals^d$: If $f(x) + f(x')$ is legal, i.e. $\{f(x),f(x')\} \neq \{\pm \infty\}$, then $f(x + x') = f(x) + f(x')$.
  \end{enumerate}
\end{definition}
If the range of $f$ is included in $\reals$, Definition \ref{def:linext} reduces to the usual definition of a linear function.
For clarity, this paper may emphasize the distinction by calling such $f$ \emph{finite linear}.

\vskip1em
To see that this definition can be satisfied nontrivially, let $f_1, f_2: \reals^d \to \reals$ be finite linear and consider $f(x) := \infty \cdot f_1(x) + f_2(x)$.
With a representation $f_1(x) = v_1 \cdot x$, we have
  \[ f(x) = \begin{cases} \infty  & v_1 \cdot x > 0  \\
                         -\infty  & v_1 \cdot x < 0  \\
                         f_2(x)   & v_1 \cdot x = 0 . \end{cases} \]

\begin{claim}
  Any such $f$ is an example of a linear extended function.
\end{claim}
\begin{myproof}
  Multiplication by a scalar distributes over the legal sum $\infty \cdot f_1(x) + f_2(x)$, so $\alpha f(x) = \infty \cdot f_1(\alpha x) + f_2(\alpha x) = f(\alpha x)$.
  To obtain additivity of $f$, consider cases on the pair $(v_1 \cdot x, ~ v_1 \cdot x')$.
  If both are zero, $f(x+x') = f_2(x+x') = f_2(x) + f_2(x') = f(x) + f(x')$.
  If they have opposite signs, then $\{f(x), f(x')\} = \{\pm \infty\}$ and the requirement is vacuous.
  If one is positive and the other nonnegative, then $v_1 \cdot (x + x') > 0$, so we have $f(x+x') = \infty = f(x) + f(x')$.
  The remaining case, one negative and the other nonpositive, is exactly analogous.
\end{myproof}

Such $f$ are not all the linear extended functions, because the subspace where $v_1 \cdot x = 0$ need not be entirely finite-valued.
As in the introduction's example, it can itself be divided into infinite and neg-infinite open halfspaces, with $f$ finite only on some further-reduced subspace.
We will show that all linear extended functions can be constructed in this recursive way, Algorithm \ref{alg:linext} (heavily relying on the setting of $\reals^d$).


\begin{algorithm}
  \caption{Computing a linear extended function $f: \reals^d \to \extreals$.}
  \label{alg:linext}
  \begin{algorithmic}
    \STATE \textbf{Parameters:} $t \in \{0,\ldots,d\}$; finite linear $\hat{f}: \reals^{d-t} \to \reals$; if $t \geq 1$, unit vectors $v_1 \in \reals^d, \ldots, v_t \in \reals^{d-t+1}$
    \STATE \textbf{Input:} $x \in \reals^d$
    \FOR{$j = 1, \ldots, t$}
      \IF{$v_j \cdot x > 0$}
        \STATE return $\infty$
      \ELSIF{$v_j \cdot x < 0$}
        \STATE return $-\infty$
      \ENDIF
      \STATE reparameterize $x$ as a vector in $\reals^{d-j}$, a member of the subspace $\{x' \st v_j \cdot x' = 0\}$
    \ENDFOR
    \STATE return $\hat{f}(x)$
  \end{algorithmic}
\end{algorithm}

First, two structural results.
Proposition \ref{prop:cvx-cnv-char}, inspired by \cite{dudik2022convex} (which uses a different construction than our axiomatic one), states that linear extended functions are exactly those that are concave, convex, and equal zero at the origin.
Lemma \ref{lemma:decomp} gives a useful decomposition of the domain of a linear extended function.
\begin{proposition}[c.f. \cite{dudik2022convex}] \label{prop:cvx-cnv-char}
  $f: \reals^d \to \extreals$ is linear extended if and only if $f(\vec{0}) = 0$ and $f$ is both concave and convex.
\end{proposition}
\begin{proof}
  $(\implies)$
  Let $f$ be linear extended.
  For any $x$, by the scaling axiom, $0 = 0 f(x) = f(0x) = f(\vec{0})$, as required.
  Now let $f(x),f(x') < \infty$, i.e. $x,x'$ are in the effective domain of $f$, and let $\rho \in [0,1]$.
  By the scaling and addition axioms, we have $\rho f(x) + (1-\rho)f(x') = f(\rho x) + f((1-\rho)x') = f(\rho x + (1-\rho)x')$, proving that $f$ is convex.
  Note that all of the above sums are legal, because none of the terms can be $+\infty$.
  Using the axioms directly, $f$ is linear extended if and only if $-f$ is linear extended, so $-f$ is also convex, so $f$ is concave as well.
  
  $(\impliedby)$
  Let $f$ be convex and concave with $f(\vec{0}) = 0$.
  We will repeatedly use that by convexity and concavity, for $\rho \in [0,1]$, we must have $f(\rho y + (1-\rho)y') = \rho f(y) + (1-\rho)f(y')$ whenever the sum is legal, which in particular holds if one summand is finite.

  We first show $f$ satisfies the scaling axiom for nonnegative coefficients.
  Let $x \in \reals^d$.
  If $\alpha = 0$, then $f(\alpha x) = f(\vec{0}) = 0 = \alpha f(x)$.
  If $\alpha \in (0,1]$, then $f(\alpha x) = f(\alpha x + (1-\alpha)\vec{0}) = \alpha f(x) + (1-\alpha)f(\vec{0}) = \alpha f(x)$, using convexity-concavity and observing that the sum is legal.
  If $\alpha > 1$, then use the previous case: $\alpha f(x) = \alpha f(\tfrac{1}{\alpha} (\alpha x)) = \alpha \tfrac{1}{\alpha}f(\alpha x) = f(\alpha x)$.

  Next, we show additivity.
  Use positive scaling and convexity-concavity to obtain $f(x + x') = 2f(\tfrac{1}{2}x + \tfrac{1}{2}x') = 2\left(\tfrac{1}{2}f(x) + \tfrac{1}{2}f(x')\right) = f(x) + f(x')$ whenever the latter sum is legal.

  Finally, we show scaling for negative coefficients, starting with $\alpha=-1$.
  By additivity, we have $0 = f(x + -x) = f(x) + f(-x)$ whenever the sum is legal.
  If the sum is legal, then $f(x)$ is finite and $f(-x) = -f(x)$.
  If the sum is illegal, then $f(x) \in \{\pm \infty\}$ and $f(-x) = -f(x)$.
  We conclude $f(-x) = -f(x)$ for all $x \in \reals^d$.
  Finally, for arbitrary $\alpha < 0$, we have $\alpha f(x) = - (-\alpha f(x)) = - f( - \alpha x) = f(- (-\alpha x)) = f(\alpha x)$, where we used positive scaling and scaling by $-1$.
\end{proof}

Now we begin to investigate the structure of linear extended functions.
Recall that a subset $S$ of $\reals^d$ is a \emph{convex cone} if, when $x,x' \in S$ and $\alpha,\beta > 0$, we have $\alpha x + \beta x' \in S$.
A convex cone need not contain $\vec{0}$.

\begin{lemma}[Decomposition] \label{lemma:decomp}
  $f: \reals^d \to \extreals$ is linear extended if and only if: (1) the sets $S^+ = f^{-1}(\infty)$ and $S^- = f^{-1}(-\infty)$ are convex cones with $S^+ = - S^-$, and (2) the set $F = f^{-1}(\reals)$ is a subspace of $\reals^d$, and (3) $f$ coincides with some finite linear function $\hat{f}$ on $F$. 
\end{lemma}
\begin{myproof}
  $(\implies)$ Let $f$ be linear extended.
  The scaling axiom implies $f(\vec{0}) = 0$, so $\vec{0} \in F$.
  Observe that $F$ is closed under scaling (by the scaling axiom) and addition (the addition axiom is never vacuous on $x,x' \in F$).
  So it is a subspace.
  $f$ satisfies scaling and additivity (never vacuous) for all members of $F$, so it is finite linear there.
  This proves (2) and (3).
  Next: the scaling axiom implies that if $x \in S^+$ then $-x \in S^-$, giving $S^+ = - S^-$ as claimed.
  Now let $x, x' \in S^+$ and $\alpha,\beta > 0$.
  The scaling axiom implies $\alpha x$ and $\beta x'$ are in $S^+$.
  The additivity axiom implies $\alpha x + \beta x' \in S^+$, proving it is a convex cone.
  We immediately get $S^- = -S^+$ is a convex cone as well, proving (1).

  $(\impliedby)$ Suppose $F$ is a subspace on which $f$ is finite linear and $S^+ = -S^-$ is a convex cone.
  We immediate obtain the useful $f(\vec{0}) = 0$ by finite linearity on a subspace (which must contain $\vec{0}$).
  We now show $f$ satisfies the two axioms needed to be linear extended.

  Scaling: if $x \in F$, the axiom follows from closure of $F$ under scaling and finite linearity of $f$ on $F$.
  Else, let $\alpha \in \reals$ and $x \in S^+$ (the case $x \in S^-$ is exactly analogous).
  If $\alpha > 0$, then $\alpha x \in S^+$ because it is a convex cone, which gives $f(\alpha x) = \alpha f(x) = \infty$.
  If $\alpha = 0$, then $f(\alpha x) = \alpha f(x) = 0$.
  If $\alpha < 0$, then use that $-x \in S^-$ by assumption, and since $S^-$ is a cone and $-\alpha > 0$, we have $(-\alpha)(-x) \in S^-$.
  In other words, $f(\alpha x) = -\infty = \alpha f(x)$, as required.

  For additivity, let $x,x'$ be given.
  \begin{itemize}
    \item If $x,x' \in F$, additivity follows from closure of $F$ under addition and finite linearity of $f$ on $F$.
    \item If $x,x' \in S^+$, then $x+x' \in S^+$ as required because it is a convex cone; analogously for $x,x' \in S^-$.
    \item If $x \in S^+, x' \in S^-$ or vice versa, the axiom is vacuously satisfied.
  \end{itemize}
  The remaining case is, without loss of generality, $x \in S^+, x' \in F$ (the proof is identical for $x \in S^-, x' \in F$).
  We must show $x + x' \in S^+$.
  Because $F$ is a subspace and $x' \in F, x \not\in F$, we must have $x + x' \not\in F$.
  Now suppose for contradiction that $x + x' \in S^-$.
  We have $-x \in S^-$ because $S^- = -S^+$.
  Because $S^-$ is a convex cone, it is closed under addition, so $x + x' + (-x) = x' \in S^-$, a contradiction.
  So $x+x' \in S^+$.
\end{myproof}

So far, our results (Proposition \ref{prop:cvx-cnv-char} and Lemma \ref{lemma:decomp}) have not actually required the domain to be finite-dimensional.
To proceed, however, we will begin to rely heavily on induction on the dimension $d$.
A useful corollary of Lemma \ref{lemma:decomp} is the following definition:
\begin{definition}[Rank] \label{def:rank}
  The \emph{rank} of a linear extended function $f: \reals^d \to \extreals$ is $d-k$, where $k$ is the dimension of the subspace $F = f^{-1}(\reals)$.
\end{definition}
In other words, the rank is the number of dimensions along which $f$ is vertical; we will show that it equals the (unique) number of vectors $t$ required to implement $f$ in Algorithm \ref{alg:linext}.
This will also imply that Definition \ref{def:rank} coincides with the definition of \emph{astral rank} in \citet{dudik2022convex}.

The effective domain of a linear extended function, i.e. the set $F \cup S^-$, at least includes an open halfspace.
If $f$ has rank zero (is finite linear), the effective domain is $\reals^d$; if rank one, it is a closed halfspace $\{ x \st v_1 \cdot x \leq 0\}$; and otherwise it is neither a closed nor open set.
The effective domain is an example of a \emph{hemispace}~\citep{martinezlegaz1988structure}, a convex set whose complement is convex.

\begin{lemma}[Recursive definition] \label{lemma:recursive}
  $f: \reals^d \to \extreals$ is linear extended if and only if one of the following hold:
  \begin{enumerate}
    \item $f$ is finite linear (this case must hold if $d=0$), or
    \item There exists a unit vector $v_1$ and linear extended function $f_2$ on the $d-1$ dimensional subspace $\{x \st v_1 \cdot x = 0\}$ such that $f(x) = f_2(x)$ if $x \cdot v_1 = 0$, else $f(x) = \infty \cdot \sign(v_1 \cdot x)$.
  \end{enumerate}
\end{lemma}
\begin{myproof}
  $(\implies)$ Suppose $f$ is linear extended.
  The case $d=0$ is immediate, as $f(\vec{0}) = 0$ by the scaling axiom.
  So let $d \geq 1$ and suppose $f$ is not finite linear; we show case (2) holds.
  Let $S^+ = f^{-1}(\infty)$, $S^- = f^{-1}(-\infty)$, and $F = f^{-1}(\reals)$.
  Recall from Lemma \ref{lemma:decomp} that $F$ is a subspace, necessarily of dimension $< d$ by assumption that $f$ is not finite; meanwhile $S^+ = -S^-$ and both are convex cones.

  We first claim that there is an open halfspace on which $f(x) = \infty$, i.e. included in $S^+$.
  First, $\cl{S^+}$ includes a closed halfspace: if not, the set $\cl{S^+} \cup \cl{S^-} \neq \reals^d$ and then its complement, an open set, would necessarily have affine dimension $d$ yet would be included in $F$, a contradiction.
  Now, because $S^+$ is convex, it includes the relative interior of $\cl{S^+}$, so it includes an open halfspace.
  
  Write this open halfspace $\{x \st v_1 \cdot x > 0\}$ for some unit vector $v_1$.
  Because $S^- = -S^+$, we have $f(x) = -\infty$ on the complement $\{x \st v_1 \cdot x < 0\}$.
  Let $f_2$ be the restriction of $f$ to the remaining set, the $d-1$ dimensional subspace $\{x \st v_1 \cdot x = 0\}$.
  Because $f$ satisfies the axioms of a linear extended function, $f_2$ is a linear extended function as well.

  $(\impliedby)$ If case (1) holds and $f$ is finite linear, then it is immediately linear extended as well.
  In case (2), will show that $f$ satisfies the conditions of Lemma \ref{lemma:decomp}.
  First, apply Lemma \ref{lemma:decomp} to the linear extended $f_2$.
  We obtain that it is finite linear on a subspace of $\{ x \st v_1 \cdot x = 0\}$, which is a subspace of $\reals^d$, giving that $f$ is finite linear on a subspace $F$.
  We also obtain $f_2^{-1}(\infty) = -f_2^{-1}(-\infty)$ and is a convex cone.
  It follows directly that $f^{-1}(\infty) = - f^{-1}(-\infty)$, i.e. $S^+ = -S^-$.
  It only remains to show that $S^+$ is a convex cone.
  In fact, $S^+ = f^{-1}(\infty) = f_2^{-1}(\infty) \cup \{ x \st v_1 \cdot x > 0\}$.
  The first set is a convex cone lying in the closure of the second set, also a convex cone.
  So the union is a convex cone: scaling is immediate; any nontrivial convex combination of a point from each set lies in the second, giving convexity; scaling and convexity give additivity.
  This shows that $S^+$ is a convex cone, the final piece needed to apply Lemma \ref{lemma:decomp} and declare $f$ linear extended.
\end{myproof}

\begin{proposition}[Correctness of Algorithm \ref{alg:linext}] \label{prop:alg-correct}
  A function $f: \reals^d \to \extreals$ is linear extended if and only if it is computed by Algorithm \ref{alg:linext} for some $t \in \{0,\ldots,d\}$, some $v_1 \in \reals^d, \ldots, v_t \in \reals^{d-t+1}$, and some finite linear $\hat{f}: \reals^{d-t} \to \reals$.
\end{proposition}
\begin{myproof}
  $(\implies)$ Suppose $f$ is linear extended.
  By Lemma \ref{lemma:recursive}, there are two cases.
  If $f$ is finite linear, then take $t=0$ and $\hat{f} = f$ in Algorithm \ref{alg:linext}.
  Otherwise, Lemma \ref{lemma:recursive} gives a unit vector $v_1$ so that $f(x) = \infty$ if $v_1 \cdot x > 0$ and $f(x) = -\infty$ if $v_1 \cdot x < 0$, as in Algorithm \ref{alg:linext}.
  $f$ is linear extended on $\{ x \st v_1 \cdot x = 0\}$, so we iterate the procedure until reaching a subspace where $f$ is finite linear, setting $t$ to be the number of iterations.

  $(\impliedby)$ Suppose $f$ is computed by Algorithm \ref{alg:linext}.
  We will use the two cases of Lemma \ref{lemma:recursive} to show $f$ is linear extended.
  If $t = 0$, then $f$ is finite linear, hence linear extended (case 1).
  If $t \geq 1$, then $f$ is in case 2 with unit vector $v_1$ and function $f_2$ equal to the implementation of Algorithm \ref{alg:linext} on $t-1$, $\hat{f}$, and $v_2,\dots,v_t$.
  This proves by induction on $t$ that, if $f$ is computed by Algorithm \ref{alg:linext}, then it satisfies one of the cases of Lemma \ref{lemma:recursive}, so it is linear extended.
\end{myproof}

\begin{figure}[ht]
  \caption{The domains of two rank-$2$ linear extended functions on $\reals^2$. In each case, the function is defined as: $f(x) = \infty \cdot \sign(v_1 \cdot x)$ if $v_1 \cdot x \neq 0$, otherwise $f(x) = \infty \cdot \sign(v_2 \cdot x)$. The shaded orange region is $f^{-1}(\infty)$, while $f(\vec{0}) = 0$ and the remaining space is $f^{-1}(-\infty)$. Coordinate axes, gray, are thickened to illustrate the behavior on those subspaces. When $f$ represents a vertical tangent plane to a convex function $g$, then $\effdom{g}$ must lie in the non-orange region.}
  \label{fig:order-matters}
  \hspace{0.15\linewidth}
  \begin{subfigure}{0.24\linewidth}
    \resizebox{\linewidth}{!}{\includegraphics{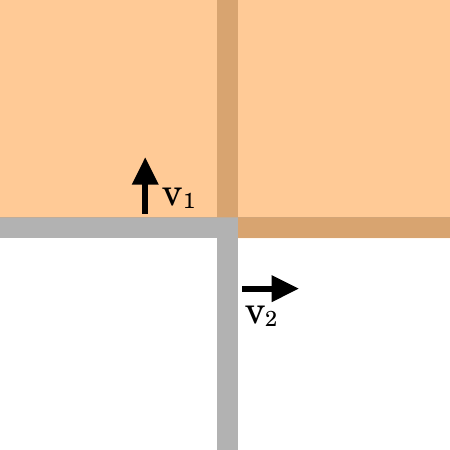}}
  \end{subfigure}
  \hfill
  \begin{subfigure}{0.24\linewidth}
    \resizebox{\linewidth}{!}{\includegraphics{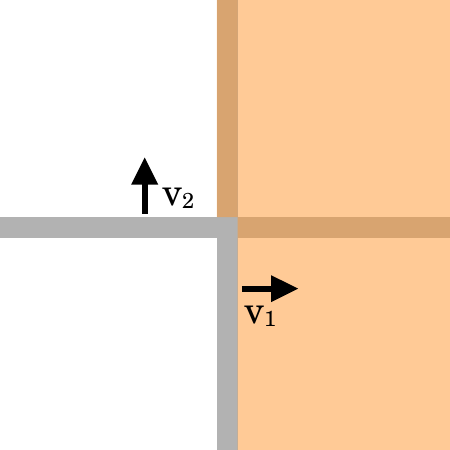}}
  \end{subfigure}
  \hspace{0.15\linewidth}
\end{figure}

The rank of $f$ is precisely its value of $t$ in Algorithm \ref{alg:linext}.
We will show next that this $t$ is unique and, in fact, the vectors $v_1,\dots,v_t$ in Algorithm \ref{alg:linext} are unique for a given $f$.
The finite-valued subspace $F = f^{-1}(\reals)$ will not change if we take these vectors in any order in Algorithm \ref{alg:linext}, reparameterizing.
However, the order determines the regions on which $f$ is infinite and neg-infinite; see Figure \ref{fig:order-matters}.
This becomes important when $f$ captures a vertical tangent plane to a convex function (Section \ref{sec:subgrad}).

\begin{proposition}[Parsimonious parameterization] \label{prop:parsim}
  Each linear extended function has a unique representation by the parameters of Algorithm \ref{alg:linext}.
\end{proposition}
\begin{myproof}
  For $i \in \{1,2\}$, let $f^{(i)}$ be the function computed by Algorithm \ref{alg:linext} with the parameters $t^{(i)}$, $\hat{f}^{(i)}$, $\{v_j^{(i)} \st j =1,\dots,t^{(i)}\}$.
  We will prove that $f^{(1)}$ and $f^{(2)}$ are distinct if any of their parameters differ: i.e. if $t^{(1)} \neq t^{(2)}$, or else $\hat{f}^{(1)} \neq \hat{f}^{(2)}$, or else there exists $j$ with $v_j^{(1)} \neq v_j^{(2)}$.

  By Lemma \ref{lemma:decomp}, each $f^{(i)}$ is finite linear on a subspace $F^{(i)}$ and positive (negative) infinite on a convex cone $S^{i,+}$ (respectively, $S^{i,-}$) with $S^{i,+} = -S^{i,-}$.
  It follows that they are equal if and only if: $S^{1,+} = S^{2,+}$ (this implies $F^{(1)} = F^{(2)}$) and they coincide on $F^{(1)}$.

  If $t^{(1)} \neq t^{(2)}$, then the dimensions of $F^{(1)}$ and $F^{(2)}$ differ, so they are nonequal, so $S^{1,+} \neq S^{2,+}$ and the functions are not the same.
  So suppose $t^{(1)} = t^{(2)}$.
  Now suppose the unit vectors are not the same, i.e. there is some smallest index $j$ such that $v^{(1)}_j \neq v^{(2)}_j$.
  Observe that on iteration $j$ of the algorithm, the $d-j+1$ dimensional subspace under consideration is identical for $f^{(1)}$ and $f^{(2)}$.
  But now there is some $x$ with $v^{(1)}_j \cdot x > 0$ while $v^{(2)}_j \cdot x < 0$, for example, $x = v^{(1)}_j - v^{(2)}_j$.
  On this $x$ (parameterized as a vector in $\reals^d$), $f^{(1)}(x) = \infty$ while $f^{(2)}(x) = -\infty$, so the functions differ.

  Finally, suppose $t^{(1)} = t^{(2)}$ and all unit vectors are the same.
  Observe that $F^{(1)} = F^{(2)}$.
  But if $\hat{f}^{(1)} \neq \hat{f}^{(2)}$, then there is a point in $F^{(1)}$ where $\hat{f}^{(1)}$ and $\hat{f}^{(2)}$ differ.
  On this point (parameterized as a vector in $\reals^d$), $f^{(1)}$ and $f^{(2)}$ differ.
\end{myproof}

While finite linear functions on $\reals^d$ are uniquely identified by $d$ real-valued parameters, linear extended functions require as many as ${d \choose 2} + 1 = \Omega(d^2)$: Unit vectors in $\reals^k$ require $k-1$ parameters (for $k \geq 2$), so even assuming rank $t=d-1$, identifying the unit vectors takes $d-1 + \cdots + 1 = {d \choose 2}$ parameters, and one more defines $\hat{f}: \reals \to \reals$.

\paragraph{Structure.}
This paper avoids topological tools, and it is not immediately clear if linear extended functions are topologically nice.
However, \citet[Ch. 4]{dudik2022convex} put a topology on these functions, coinciding with pointwise convergence, in which they form a compactification of the set of linear functions.
That work also considers algebraic structure, defining a leftward addition operator and notions of convex sets in the space of linear extended functions (Ch. 9).

\section{Extended subgradients} \label{sec:subgrad}

Recall that a \emph{subgradient} of a function $g$ at a point $x_0$ is a (finite) linear function $f$ satisfying $g(x) \geq g(x_0) + f(x-x_0)$ for all $x$.
The following definition of extended subgradient replaces $f$ with a linear extended function; to force sums to be legal, we will only define it at points where the function is finite.
A main aim of this section is to show that a function, having a convex effective domain, is convex if and only if it has an extended subgradient everywhere it is finite (Proposition \ref{prop:subgrad-char}).

\begin{definition}[Extended subgradient] \label{def:extsub}
  Given a function $g: \reals^d \to \extreals$, a linear extended function $f$ is an \emph{extended subgradient} of $g$ at a point $x_0$ where $g(x_0) \in \reals$ if, for all $x \in \reals^d$, $g(x) \geq g(x_0) + f(x - x_0)$.
\end{definition}
Again for clarity, if this $f$ is a finite linear function, we may call it a \emph{finite} subgradient of $g$.
Note that the sum appearing in Definition \ref{def:extsub} is always legal because $g(x_0) \in \reals$.

\vskip1em
As is well-known, finite subgradients correspond to supporting hyperplanes of the epigraphs of convex functions $g$, but they do not exist at points where all such hyperplanes are vertical; consider, in one dimension,
\begin{equation}
    g(x) = \begin{cases}
            0  &  x < 1  \\
            1  &  x = 1  \\
            \infty &  x > 1.  \end{cases} \label{eqn:example-1d-g}
\end{equation}
We next show that every convex $g$ has an \emph{extended} subgradient everywhere it is finite.
For example, $f(z) = \infty \cdot z$ is an extended subgradient of the above $g$ at $x=1$.
\begin{proposition}[Existence of extended subgradients] \label{prop:subgrad-exist}
  Each convex function $g: \reals^d \to \extreals$ has an extended subgradient at every point $x$ where $g(x) \in \reals$.
\end{proposition}
\begin{myproof}
  Let $g(x) \in \reals$; we construct an extended subgradient $f_x$.
  Define $g_x(z) = g(x+z) - g(x)$, the shift of $g$ that moves $(x,g(x))$ to $(\vec{0},0)$.
  Note $g_x$ is convex, as its epigraph is a shift of $g$'s.
  Now we appeal to a technical lemma, Lemma \ref{lemma:conv-0-minorized} (next), which says there exists a linear extended $f_x$ that minorizes $g_x$.
  This completes the proof: For any $x' \in \reals^d$, we have $f_x(x'-x) \leq g_x(x'-x)$, which rearranges to give $g(x) + f_x(x'-x) \leq g(x')$, with legality from $g(x) \in \reals$.
\end{myproof}

\begin{lemma} \label{lemma:conv-0-minorized}
  Let $g_x : \reals^d \to \extreals$ be a convex function with $g_x(\vec{0}) = 0$.
  Then there exists a linear extended function $f_x$ with $f_x(z) \leq g_x(z)$ for all $z \in \reals^d$.
\end{lemma}
\begin{myproof}
  First we make Claim (*): if $\vec{0}$ is in the relative interior of the effective domain $G_x$ of $g_x$, then such an $f_x$ exists.
  This follows because, in this case, $g_x$ is a proper function (an improper function is $-\infty$ in the relative interior of its domain \citep[Theorem 7.2]{rockafellar1970convex}); and the proper function $g_x$ necessarily has at relative interior points a finite subgradient \citep[Theorem 23.4]{rockafellar1970convex}, which we can take to be $f_x$.

  Now, we prove the result by induction on $d$.

  If $d=0$, then the only linear extended function is $f_x = g_x$.

  If $d \geq 1$, there are two cases.
  If $\vec{0}$ is in the interior of $G_x$, then Claim (*) applies and we are done.

  Otherwise, $\vec{0}$ is in the relative boundary of the convex set $G_x = \effdom{g_x}$.
  This implies\footnote{Definition and existence of a supporting hyperplane are referred to e.g. in \citet[A]{hiriarturrut2001fundamentals}: Definition 1.1.2 of a hyperplane (there $v$ must merely be nonzero, but it is equivalent to requiring a unit vector); Definition 2.4.1 of a supporting hyperplane to $G_x$ at $x$, specialized here to the case $x = \vec{0}$; Lemma 4.2.1, giving existence of a supporting hyperplane to any nonempty convex set $G_x$ at any boundary point, which includes the relative boundary.} that there exists a hyperplane supporting $G_x$ at $\vec{0}$.
  That is, there is some unit vector $v \in \reals^d$ such that $G_x \subseteq \{z \st v \cdot z \leq 0\}$.
  Set $f_x(z) = \infty$ if $v \cdot z > 0$ and $f_x(z) = -\infty$ if $v \cdot z < 0$.
  Note $f_x$ minorizes $g$ on these regions, in the first case because $g_x(z) = \infty = f_x(z)$ (by definition of effective domain), and in the second case because $f_x(z) = -\infty$.

  So it remains to define $f_x$ on the subspace $S := \{ z \st v \cdot z= 0\}$ and show that it minorizes $g_x$ there.
  But $g_x$, restricted to this subspace, is again a proper convex function equalling $0$ at $\vec{0}$, so by induction we have a minorizing linear extended function $\hat{f_x}$ on this $d-1$ dimensional subspace.
  Define $f_x(z) = \hat{f_x}(z)$ for $z \in S$.
  Then $f_x$ minorizes $g_x$ everywhere.
  We also have that $f_x$ is linear extended by Lemma \ref{lemma:recursive}, as it satisfies the necessary recursive format.
\end{myproof}
This proof would have gone through if Claim (*) used ``interior'' rather than ``relative interior'', because supporting hyperplanes exist at the entire boundary, not just the relative boundary.
That version would have constructed extended subgradients of possibly greater rank (Definition \ref{def:rank}).

Examples of rank-$2$ subgradients $f$ constructed by Lemma \ref{lemma:conv-0-minorized} appeared in Figure \ref{fig:order-matters}.
Concretely, in those examples we can suppose $g: \reals^2 \to \extreals$ equals $0$ at the origin and $-1$ elsewhere in its effective domain, with $f$ representing the subgradient at the origin.
Here, $g$'s epigraph has a vertical tangent plane above the subspace $\{x : v_1 \cdot x = 0\}$.
Then $g$ restricted to that subspace is a convex one-dimensional function; and this function's epigraph has a vertical tangent line at the origin.

\vskip1em
It is now useful to define \emph{affine extended functions}.
Observe that a bad definition would be ``a linear extended function plus a constant.''
Both vertical and \emph{horizontal} shifts of linear extended functions must be allowed in order to capture, e.g.
\begin{equation}
     h(x) = \begin{cases} -\infty & x < 1  \\
                          1       & x = 1  \\
                          \infty  & x > 1 . \end{cases} \label{eqn:example-1d-h}
\end{equation}

\begin{definition}[Affine extended function] \label{def:affext}
  A function $h: \reals^d \to \extreals$ is \emph{affine extended} if for some $\beta \in \reals$ and $x_0 \in \reals^d$ and some linear extended $f: \reals^d \to \extreals$ we have $h(x) = f(x - x_0) + \beta$.
\end{definition}
\begin{definition}[Supports] \label{def:supports}
  An affine extended function $h$ \emph{supports} a function $g: \reals^d \to \extreals$ at a point $x_0 \in \reals^d$ if $h(x_0) = g(x_0)$ and $h$ minorizes $g$, i.e. $h(x) \leq g(x)$ for all $x \in \reals^d$.
\end{definition}
For example, the $h$ in Display \ref{eqn:example-1d-h} supports, at $x_0=1$, the convex function $g$ of Display \ref{eqn:example-1d-g}.
Actually, $h$ supports $g$ at all $x_0 \geq 1$.

\begin{observation} \label{obs:affine}
  Given a linear extended $f$, some $\beta \in \reals$, and some $x_0 \in \reals^d$, define the affine extended function $h(x) = f(x-x_0) + \beta$.
  \begin{enumerate}
    \item $h$ is convex and concave.
    \item Let $g: \reals^d \to \extreals$ such that $g(x_0) \in \reals$.
          Then $h$ supports $g$ at $x_0$ if and only if $f$ is an extended subgradient of $g$ at $x_0$.
    \item If $h(x_1)$ is finite, then for all $x$ we have $h(x) = f(x-x_1) + \beta_1$ where $\beta_1 = f(x_1-x_0) + \beta$.
  \end{enumerate}
\end{observation}
Point 1 follows because $f$ is convex (respectively, concave), and $h$'s epigraph (respectively, hypograph) is a shift of $f$'s.
Point 2 follows because $h$ supports $g$ at $x_0$ if and only if $\beta = g(x_0)$ and $g(x) \geq \beta + f(x-x_0)$ for all $x$.
Point 3 follows because $h(x_1) = f(x_1 - x_0) + \beta$, so $f(x_1 - x_0)$ is finite: the sum $f(x-x_1) + f(x_1-x_0)$ is always legal and equals $f(x-x_0)$.

\vskip1em
An important fact in convex analysis is that a \emph{closed} (i.e. lower semicontinuous) convex function $g$ is the pointwise supremum of a family of affine functions.
Restating, closed convex epigraphs are the intersection of a family of non-vertical closed halfspaces.
Using affine extended functions, we can drop the \emph{non-vertical} requirement and recover all convex epigraphs, i.e. characterize all convex functions.
(They are no longer closed because a vertical tangent plane does not enclose a closed set, i.e. the epigraph of a linear extended function of nonzero rank is not closed.)

\begin{proposition}[Supremum characterization] \label{prop:sup-char}
  A function $g: \reals^d \to \extreals$ is convex if and only if it is the pointwise supremum of a family of affine extended functions.
\end{proposition}
\begin{myproof}
  $(\impliedby)$
  Let $H$ be a family of affine extended functions and let $g(x) := \sup_{h \in H} h(x)$.
  The epigraph of each affine extended function is a convex set (Observation \ref{obs:affine}), and $g$'s epigraph is the intersection of these, hence a convex set, so $g$ is convex.

  $(\implies)$ Suppose $g$ is convex; let $H$ be the set of affine extended functions minorizing $g$.
  We show $g(x) = \sup_{h \in H} h(x)$.

  First consider $x$ in the effective domain $G$ of $g$.
  Suppose $g(x) = -\infty$.
  If $H = \emptyset$, then $g(x) = \sup \emptyset = -\infty$ as required.
  Otherwise, every minorizing $h$ has $h(x) = -\infty$, so $-\infty = \sup_{h \in H} h(x)$.
  Now suppose $g(x) \in \reals$.
  By Proposition \ref{prop:subgrad-exist}, $g$ has an extended subgradient at $x$, which implies (Observation \ref{obs:affine}) it has a supporting affine function $h_x$ there, and of course $h_x \in H$.
  We have $g(x) = h_x(x)$ and by definition $g(x) \geq h(x)$ for all $h \in H$, so $g(x) = \max_{h \in H} h(x)$.

  Second, consider $x \not\in G = \effdom{g}$.
  We must show that $\sup_{h \in H} h(x) = \infty$.
  We apply a technical lemma, Lemma \ref{lemma:sup-affine-indicates} (stated and proven next), to obtain a set $H'$ of affine extended functions that all equal $-\infty$ on $G$ (hence $H' \subseteq H$), but for which $\sup_{h \in H'} h(x) = \infty$, as required.
\end{myproof}

\begin{lemma} \label{lemma:sup-affine-indicates}
  Let $G$ be a convex set in $\reals^d$ and let $x \not\in G$.
  There is a set $H'$ of affine extended functions such that $\sup_{h \in H'} h(x) = \infty$ while $h(x') = -\infty$ for all $h \in H$ and $x' \in G$.
\end{lemma}
For intuition, observe that the easiest case is if $x \not\in \cl{G}$, when we can use a strongly separating hyperplane to get a single $h$ that is $\infty$ at $x$ and $-\infty$ on $G$.
The most difficult case is when $x$ is an extreme, but not exposed, point of $\cl{G}$.
(Picture a convex $g$ with effective domain $\cl{G}$; now modify $g$ so that $g(x) = \infty$.)
To capture this case, we use the full recursive structure of linear extended functions.
\begin{myproof}
  Let $G' = \{x' - x \st x' \in G\}$.
  Note that $G'$ is a convex set not containing $\vec{0}$.
  By Lemma \ref{lemma:linext-indicates}, next, there is a linear extended function $f$ equalling $-\infty$ on $G'$.
  Construct $H = \{ h(x') = f(x'-x) + \beta \st \beta \in \reals\}$.
  Then $\sup_{h \in H} h(x) = \sup_{\beta \in \reals} \beta = \infty$, while for $x' \in G$ we have $h(x') = f(x'-x) + \beta = -\infty$.
  This proves Lemma \ref{lemma:sup-affine-indicates}.
\end{myproof}

\begin{lemma} \label{lemma:linext-indicates}
  Given a convex set $G' \subseteq \reals^d$ not containing $\vec{0}$, there exists a linear extended $f$ with $f(x') = -\infty$ for all $x' \in G'$.
\end{lemma}
\begin{myproof}
  Note that if $G' = \emptyset$, then the claim is immediately true; and if $d=0$, then $G'$ must be empty.
  Now by induction on $d$:
  If $d=1$, then we can choose $f(z) = \infty \cdot z$ or else $f(z) = -\infty \cdot z$; one of these must be $-\infty$ on the convex set $G'$.

  Now suppose $d \geq 2$; we construct $f$.
  The convex hull of $G' \cup \{\vec{0}\}$ must be supported at $\vec{0}$, which is a boundary point\footnote{If $\vec{0}$ were in the interior of this convex hull, then it would be a convex combination of other points in the convex hull, so it would be a convex combination of points in $G'$, a contradiction. For the supporting hyperplane argument, see the above footnote, using reference \citet{hiriarturrut2001fundamentals}.}, by some hyperplane.
  Let this hyperplane be parameterized by the unit vector $v_1$ where $G' \subseteq \{z \st v_1 \cdot z \leq 0\}$.
  Set $f(z) = \infty$ if $v_1 \cdot z > 0$ and $-\infty$ if $v_1 \cdot z < 0$.
  It remains to define $f$ on the subspace $S := \{z \st v_1 \cdot z = 0\}$.
  Immediately, $G' \cap S$ is a convex set not containing $\vec{0}$.
  So by inductive hypothesis, there is a linear extended $\hat{f}: \reals^{d-1} \to \extreals$ that is $-\infty$ everywhere on $G' \cap S$.
  Set $f(z) = \hat{f}(z)$ on $S$.
  Now we have $f$ is linear extended because it satisfies the recursive hypotheses of Lemma \ref{lemma:recursive}.
  We also have that $f$ is $-\infty$ on $G'$: each $z \in G'$ either has $v_1 \cdot z < 0$ or $z \in S$, and both cases have been covered.
\end{myproof}
Consider the convex function in one dimension,
\begin{align}
  g(x) = \begin{cases}
         0  &  x < 1   \\
         \infty  &  x \geq 1 . \end{cases} \label{eqn:example-open}
\end{align}
For the previous proof to obtain $g$ as a pointwise supremum of affine extended functions, it needed a sequence of the form
  \[  h(x) = \begin{cases}
                -\infty  &  x < 1   \\
                \beta    &  x = 1   \\
                \infty   &  x > 1  \end{cases}  \]
for $\beta \to \infty$.
Indeed, one can show that $g$ has no affine extended supporting function at $x=1$.
But it does everywhere else: in particular, any such $h$ supports $g$ at all $x > 1$.
So we might have hoped that each convex function is the pointwise \emph{maximum} of affine extended functions; but the example of Display \ref{eqn:example-open} and $x=1$ prevents this.
This hope does hold on $g$'s finite-valued domain, implied by existence of extended subgradients there (Proposition \ref{prop:subgrad-exist}).
Speaking of which, we are now ready to prove a converse.
\begin{proposition}[Subgradient characterization] \label{prop:subgrad-char}
  A function $g: \reals^d \to \extreals$ is convex if and only if (a) its effective domain is convex and (b) it has an extended subgradient at every point $x$ where $g(x) \in \reals$.
\end{proposition}
\begin{myproof}
  $(\implies)$
  This is Proposition \ref{prop:subgrad-exist}, along with convexity of $\effdom{g}$ for convex $g$.

  $(\impliedby)$
  Suppose $G = \effdom{g}$ is convex and $g$ has an extended subgradient at each point where $g(x) \in \reals$.
  We will show $g$ is the pointwise supremum of a family $H$ of affine extended functions, hence (Proposition \ref{prop:sup-char}) convex.

  First, let $H^*$ be the set of affine extended functions that support $g$.
  Then $H^*$ must include a supporting function at any $x$ with $g(x) \in \reals$ (Observation \ref{obs:affine}).
  Now we repeat a trick: for each $x \not\in G$, we obtain from Lemma \ref{lemma:sup-affine-indicates} a set $H_x$ of affine extended functions that are all $-\infty$ on $G$, where $\sup_{h \in H_x} h(x) = \infty$.
  Letting $H = H^*$ unioned with $\bigcup_{x \not\in G} H_x$, we claim $g(x) = \sup_{h \in H} h(x)$.
  First, every $h \in H$ minorizes $g$: this follows immediately for $h \in H^*$, and by definition of effective domain for each $H_x$.
  So $g(x) \geq \sup_{h \in H} h(x)$.
  If $g(x) = \infty$, equality follows by construction of $H_x$.
  If $g(x) \in \reals$, equality follows because $H^*$ contains a supporting function with $h(x) = g(x)$.
  If $g(x) = -\infty$, equality follows because for all $h \in H$, $-\infty = g(x) \geq h(x) = -\infty$.
\end{myproof}
It should be possible to remove the requirement of a convex effective domain in Proposition \ref{prop:subgrad-char}, if we broaden the definition of extended subgradients to exist at points where $g(x) \in \{-\infty,\infty\}$.
This is the case in \cite[Ch 18.2]{dudik2022convex}.
Under such a formalization, convex functions would be characterized by having extended subgradients everywhere on $\reals^d$.
However, this would break the bijection (Observation \ref{obs:affine} (2)) between extended subgradients and supporting extended affine functions, which may not exist at all points -- e.g., the example of Display \ref{eqn:example-open}.

\vskip1em
Finally, characterizations of strict convexity will be useful for scoring rule applications.
We consider proper convex functions, where there exist extended subgradients and affine extended supports everywhere on $\effdom{g} = g^{-1}(\reals)$.
Say an affine extended function is \emph{uniquely supporting} of $g$ at $b \in \P \subseteq \effdom{g}$, with respect to $\P$, if it supports $g$ at $b$ and at no other point in $\P$.

\begin{lemma}[Strict convexity] \label{lemma:strict-subgrad}
  Let $g: \reals^d \to \reals \cup \{\infty\}$ be convex, and let $\P \subseteq \effdom{g}$ be a convex set.
  The following are equivalent: (1) $g$ is strictly convex on $\P$; (2) no two distinct points in $\P$ share an extended subgradient; (3) every affine extended function supports $g$ in $\P$ at most once; (4) $g$ has at each $x \in \P$ a uniquely supporting affine extended function.
\end{lemma}
Although the differences between (2,3,4) are small, they may be useful in different scenarios.
For example, proving (4) may be the easiest way to show $g$ is strictly convex, whereas (3) may be the more powerful implication of strict convexity.
\begin{myproof}
  We prove a ring of contrapositives.

  $(\lnot 1 ~ \implies ~ \lnot 4)$
  If $g$ is convex but not strictly convex on $\P \subseteq \effdom{g}$, then there are two points $a,c \in \P$ with $a \neq c$ and there is $0 < \rho < 1$ such that, for $b = \rho \cdot a + (1-\rho)c$, $g(b) = \rho g(a) + (1-\rho)) g(c)$.
  We show $g$ has no uniquely supporting $h$ at $b$; intuitively, tangents at $b$ must be tangent at $a$ and $c$ as well.

  Let $h$ be any supporting affine extended function at $b$ and write $h(x) = g(b) + f(x-b)$.
  Because $h$ is supporting, $g(c) \geq h(c) = g(b) + f(c-b)$, implying with the axioms of linear extended functions that $f(b-c) \geq g(b) - g(c)$.
  Now, $b-c = \rho(a-c) = \frac{\rho}{1-\rho}(a-b)$.
  So $f(b-c) = \frac{\rho}{1-\rho}f(a-b)$.
  Similarly, $g(b) - g(c) = \rho(g(a) - g(c)) = \frac{\rho}{1-\rho}(g(a) - g(b))$.
  So $f(a-b) \geq g(a) - g(b)$, so $h(a) = g(b) + f(a-b) \geq g(a)$.
  But because $h$ supports $g$, we also have $h(a) \leq g(a)$, so $h(a) = g(a)$ and $h$ supports $g$ at $a$.

  $(\lnot 4 ~ \implies ~ \lnot 3)$
  Almost immediate: $g$ has a supporting affine extended function at every point in $\P$ because it has an extended subgradient everywhere (Proposition \ref{prop:subgrad-exist}).
  If $\lnot 4$, then one supports $g$ at two points in $\P$.

  $(\lnot 3 ~ \implies ~ \lnot 2)$
  Suppose an affine extended $h$ supports $g$ at two distinct points $a,b \in \P$.
  Using Observation \ref{obs:affine}, we can write $h(x) = f(x - a) + g(a) = f(x-b) + g(b)$ for a unique linear extended $f$, so $f$ is an extended subgradient at two distinct points.

  $(\lnot 2 ~ \implies ~ \lnot 1)$
  Suppose $f$ is an extended subgradient at distinct points $a,c \in \P$.
  By definition of extended subgradient, we have $g(a) \geq g(c) + f(a-c)$ and $g(c) \geq g(a) - f(a-c)$, implying $f(a-c) = g(a) - g(c)$.

  Now let $b = \rho \cdot a + (1-\rho) c$ for any $0 < \rho < 1$.
  By convexity and linear extended axioms, $g(b) \geq g(a) - f(a-b) = g(a) - f((1-\rho)(a-c)) = g(a) - (1-\rho)(g(a) - g(c)) = \rho g(a) + (1-\rho)g(c)$.
  So $g$ is not strictly convex. 
\end{myproof}

\paragraph{The extended subdifferential.}
Although its topological properties and general usefulness are unclear, we should probably not conclude without defining the \emph{extended subdifferential} of the function $g: \reals^d \to \extreals$ at $x$ to be the set of extended subgradients of $g$ at $x$, nonempty if $g(x) \in \reals$.
On the interior of $\effdom{g}$, the extended subdifferential can only contain finite subgradients.
But note that if the effective domain has affine dimension smaller than $\reals^d$, then the epigraph is contained in one of its vertical tangent supporting hyperplanes.
In this case, the relative interior, although it has finite subgradients, will have extended ones as well that take on non-finite values outside the effective domain.

%

\section{Proper scoring rules} \label{sec:proper-scoring-rules}

This section will first define scoring rules; then recall their characterization and discuss why it is not fully constructive; then use linear extended functions to prove more complete and constructive versions.
In particular, given any subset $\P$ of a finite-dimensional probability simplex, Theorem \ref{thm:proper-subset-char} constructively characterizes the convex functions (or \emph{entropies}) that give rise to a proper scoring rule.

\subsection{Definitions}
$\Y$ is a finite set of mutually exclusive and exhaustive outcomes, also called observations.
For a vector $v \in \reals^{\Y}$, we write $v(y)$ for the $y$th component of $v$.
The probability simplex on $\Y$ is $\Delta_{\Y} = \{p \in \reals^{\Y} \st (\forall y) ~ p(y) \geq 0 ~,~ \sum_{y \in \Y} p(y) = 1\}$.
The \emph{support} of $p \in \Delta_{\Y}$ is $\Supp{p} = \{y : p(y) > 0\}$.
The distribution with $p(y)=1$ is denoted $\delta_y \in \Delta_{\Y}$.

Proper scoring rules model an expert providing a forecast $p \in \Delta_{\Y}$, after which $y \in \Y$ is observed and the expert's score is $S(p,y)$.
The expert chooses $p$ to maximize expected score according to an internal belief $q \in \Delta_{\Y}$.
More generally, it is sometimes assumed that reports and beliefs are restricted to a subset $\P \subseteq \Delta_{\Y}$.
Two of the most well-known are the log scoring rule, where $S(p,y) = \log p(y)$; and the quadratic scoring rule, where $S(p,y) = -\|\delta_y - p\|_2^2$.

\begin{definition}[Scoring rule, regular] \label{def:scoring-rule}
  Let $\Y$ be finite and $\P \subseteq \Delta_{\Y}$, nonempty.
  A function $S: \P \times \Y \to \extreals$ is termed a \emph{scoring rule}.
  It is \emph{regular} if $S(p,y) \neq \infty$ for all $p \in \P$, $y \in \Y$.
\end{definition}

Regular scoring rules are nice for two reasons.
First, calculating expected scores such as $\sum_y q(y) S(p,y)$ may lead to illegal sums if $S(p,y) \in \extreals$, but regular scores guarantee that the sum is legal (using that $q(y) \geq 0$).
Second, even if illegal sums were no problem, allowing scores of $S(p,y) = \infty$ leads to strange and unexciting rules: $p$ is an optimal report for any belief $q$ unless $q(y) = 0$.
In contrast, allowing $S(p,y) = -\infty$ in some cases can be reasonable; for example, this is the case for the popular log score when $p(y) = 0$, but $y$ is observed.

\begin{definition}[Expected score, strictly proper] \label{def:proper}
  Given a regular scoring rule $S: \P \times \Y \to \reals \cup \{-\infty\}$, the \emph{expected score} for report $p \in \P$ under belief $q \in \Delta_{\Y}$ is written $S(p;q) := \sum_{y \in \Y} q(y) S(p,y)$.
  The regular scoring rule $S$ is \emph{proper} if for all $p,q \in \P$ with $p \neq q$, $S(p;q) \leq S(q;q)$, and it is \emph{strictly proper} if this inequality is always strict.
\end{definition}
The log and quadratic scoring rules mentioned above are strictly proper.
An example of a scoring rule that is not proper is $S(p,y) = p(y)$.
For more background, we refer the reader to \citet{gneiting2007strictly}.

\vskip1em
Why use a general definition of scoring rules if the only interesting rules are regular?
It will be useful for technical reasons: We will consider ``subtangent'' constructions that always yield a well-defined scoring rule, which is proper if and only if it is regular.
We can then focus on conditions under which the scoring rule is regular.
Understanding these conditions is a main contribution of Theorem \ref{thm:proper-char} and particularly Theorem \ref{thm:proper-subset-char}.

\vskip1em
\paragraph{Subtangent rules and entropies.}
The well-known proper scoring rule characterization, discussed next, roughly states that all regular proper scoring rules can be written in the following form.
\begin{definition}[Subtangent rule] \label{def:subtangent-rule}
  If a scoring rule $S: \P \times \Y \to \extreals$ satisfies
  \begin{align}
    S(p,y) &= g(p) + f_p(\delta_y - p)   & \forall p \in \P, \forall y \in \Y \label{eqn:subtan}
  \end{align}
  for some convex function $g: \reals^{\Y} \to \reals \cup \{\infty\}$ with $\P \subseteq \effdom{g}$ and some selection of its extended subgradients, $f_p$ at each $p \in \P$, then call $S$ a \emph{subtangent rule (of $g$)}.
\end{definition}
The geometric intuition is that the affine extended function $q \mapsto g(p) + f_p(q - p)$ is a linear approximation of $g$ at the point $p$.
So a subtangent rule, on prediction $p$ and observation $y$, evaluates this linear approximation at $\delta_y$.

\begin{observation} \label{obs:entropy}
  If $S$ is a regular subtangent rule of $g$, then its expected score function on $\P$ is $S(p;q) = g(p) + f_p(q-p)$, and in particular the expected score for truthful reporting is $S(p;p) = g(p)$.
\end{observation}
The observation follows because, by regularity, the expected score summation is legal, and we can use the scaling and additivity axioms to move the expectation inside $f_p$.
This implies that $g(p)$ can be viewed as a \textbf{measure of value of information} or a \textbf{generalized entropy}, where $g(p)$ captures the expected score or amount of knowledge of someone with belief $p$~\citep{howard1966information,grunwald2004game}.
For example, with the log scoring rule $S(p,y) = \log p(y)$, the corresponding $g$ is precisely the negation of Shannon entropy.

\paragraph{Subgradient rules.}
There is a less-often used, but tighter version of the characterization~\citep{mccarthy1956measures,hendrickson1971proper}: all regular proper scoring rules can be written as \emph{subgradient rules} of a \emph{positively homogeneous}\footnote{$g$ is positively homogeneous if $g(\alpha x) = \alpha g(x)$ for all $\alpha \geq 0, x \in \reals^d$.} convex function $g$, where:
\begin{definition}[Subgradient scoring rule] \label{def:subgrad-score}
  If a scoring rule $S: \P \times \Y \to \extreals$ satisfies
  \begin{align}
    S(p,y) &= f_p(\delta_y)   &  \forall p \in \P, y \in \Y  \label{eqn:subgrad-score}
  \end{align}
  for some positively homogeneous convex function $g: \reals^{\Y} \to \reals \cup \{\infty\}$ with $\P \subseteq \effdom{g}$ and some selection of its extended subgradients, $f_p$ at each $p \in \P$, then call $S$ a \emph{subgradient rule (of $g$)}.
\end{definition}

\begin{observation} \label{obs:subgradient-subtangent}
  Every regular subgradient rule $S(p,y) = f_p(\delta_y)$ of a positively homogeneous $g$ is a subtangent rule of $g$.
\end{observation}
\begin{proof}
  We must show $f_p(\delta_y) = g(p) + f_p(\delta_y - p)$ for all $p \in \P, y \in \Y$.
  We will first show that $g(p) = f_p(p)$.
  Using the definition of subgradient and homogeneity, $g(p) + f_p(p) = g(p) + f_p(2p-p) \leq g(2p) = 2g(p)$, and since $g(p) \in \reals$, we conclude $f_p(p) \leq g(p)$.
  Similarly, $g(p) - f_p(p) = g(p) + f_p(\vec{0} - p) \leq g(\vec{0}) = 0$, so $g(p) \leq f_p(p)$.
  Since $g(p) = f_p(p) \in \reals$, we have $g(p) + f_p(\delta_y - p) = g(p) + f_p(\delta_y) - f_p(p) = f_p(\delta_y)$.
\end{proof}

A \emph{characterization} theorem is therefore strongest when it concludes that any proper scoring rule can be written as a subgradient rule; this implies that it is also a subtangent rule.
On the other hand, \emph{constructions} of proper scoring rules, Theorems \ref{thm:proper-char} and \ref{thm:proper-subset-char}, are strongest when they are able to construct, not only subgradient rules from homogeneous convex functions, but also subtangent rules from general convex functions.
In other words, it is often more practically convenient to construct subtangent rules.

The distinction between subgradient and subtangent rules is not overly significant: given any convex $g$ and subtangent rule $S$ of $g$, we can define the homogeneous $g'$ to match $g$ on the simplex, and scale $g$ appropriately outside the simplex, i.e. $g'(x) = \|x\|_1 g(\tfrac{x}{\|x\|_1})$, and $g'(\vec{0}) = 0$, where $\|x\|_1 = \sum_y |x(y)|$.
Then $S$ can be expressed as a subgradient rule of $g'$.

However, the distinction has historically caused confusion.
For example, the log scoring rule $S(p,y) = \log p(y)$ is typically presented as a \emph{subtangent} rule of the negative entropy $g(p) = \sum_{y \in \Y} p(y) \log p(y)$ (infinite if $p \not\in \Delta_{\Y}$), but it is not a \emph{subgradient} rule of that function, as pointed out by \citet{marschak1959remarks}.
But, as \citet{hendrickson1971proper} responds, the log score is a subgradient rule of the positively homogeneous function $g(p) = \sum_{y \in \Y} p(y) \log \left( \tfrac{p(y)}{\|p\|_1} \right)$.

\subsection{Prior characterizations and missing pieces}

The scoring rule characterization has appeared in many forms, notably in \citet{mccarthy1956measures, savage1971elicitation, hendrickson1971proper, schervish1989general, gneiting2007strictly}.
In particular, \citet{savage1971elicitation} proved:
\begin{theorem}[\citet{savage1971elicitation} characterization] \label{thm:savage-char}
  A scoring rule $S$ on $\P = \Delta_{\Y}$ \emph{taking only finite values} is proper (respectively, strictly proper) if and only if it is of the form (\ref{eqn:subtan}) for some convex (respectively, strictly convex) $g$ with \emph{finite linear} subgradients $\{f_p : p \in \P\}$.
\end{theorem}

There are several nonconstructive wrinkles in this result, and all other characterizations of which I am aware.
By itself, the theorem above does not directly answer the following questions (and the fact that it does not may surprise even experts).
\begin{enumerate}
  \item \emph{(Tightness of strict convexity.)} Let the finite-valued $S$ be strictly proper.
        Theorem \ref{thm:savage-char} asserts that it is of the form (\ref{eqn:subtan}) for some strictly convex $g$.
        Can it also be of the form (\ref{eqn:subtan}) for some other, \emph{non}-strictly convex $g$?
  \item \emph{(Construction)} Let $g: \Delta_{\Y} \to \reals$ be convex.
        Theorem \ref{thm:savage-char} implies that \emph{if} there exists a finite-valued scoring rule $S$ of the form (\ref{eqn:subtan}), \emph{then} $S$ is proper.
        But does such an $S$ necessarily exist?
  \item \emph{(Construction on general $\P$)} What if the same question is asked, but for a domain $\P \subseteq \Delta_{\Y}$?
\end{enumerate}
For the finite-valued case, answers are not too difficult to derive, although I do not know of a citation.\footnote{%
    The answer to the first question is \emph{no}.
    This follows because a convex subdifferentiable $g$ is strictly convex on $\Delta_{\Y}$ if and only if each subgradient appears at most at one point $p \in \Delta_{\Y}$; this is an easier version of Lemma \ref{lemma:strict-subgrad}.
    The second question is deeper, and the answer is \emph{no, not necessarily}.
    Some convex $g$ are not subdifferentiable, i.e. they have no finite subgradient at some points, so it is not possible to construct a finite-valued scoring rule from them.
    We will see that the answer is \emph{yes} if we generalize Theorem \ref{thm:savage-char} to allow $S$ to be sometimes $-\infty$.
    The third question is deeper yet.
    The answer is \emph{not necessarily} even in the general case, and Theorem \ref{thm:proper-subset-char} will characterize the cases where such $S$ exist.}
But as discussed next, these questions are also unanswered by characterizations that allow scores of $-\infty$; although to an extent the answers may be known to experts.
Theorems \ref{thm:proper-char} and \ref{thm:proper-subset-char} will address them, for the case of finite $\Y$.

\paragraph{Extensions to scores with $-\infty$.}
\citet{savage1971elicitation} disallows (intentionally, Section 9.4) the log scoring rule of \citet{good1952rational}, $S(p,y) = \log p(y)$, because it can assign score $-\infty$ if $p(y) = 0$.
The prominence of the log scoring rule has historically (e.g. \citet{hendrickson1971proper}) motivated including it, but doing so requires generalizing the characterization to include possibly-neg-infinite scores.

The modern treatment of \citet{gneiting2007strictly} (also \citet{frongillo2014general,frongillo2021general}) captures such scoring rules as follows: It briefly defines extended-valued \emph{subtangents}, analogues of this paper's extended subgradients (but possibly on infinite-dimensional spaces, so for instance requirements of ``legal sums'' are replaced with ``quasi-integrable'').
The characterization can then be stated: A regular scoring rule is (strictly) proper if and only if it is of the form (\ref{eqn:subtan}) for some (strictly) convex $g$ with \emph{subtangents} $\{f_p : p \in \P\}$.

However, these characterizations have drawbacks analogous to those enumerated above, and also utilize the somewhat-mysterious subtangent objects.
This motivates a more explicit and constructive approach of this paper, albeit for finite-dimensional spaces only.

\paragraph{Contributions.}
We will see two results on the \emph{construction of proper scoring rules}, extensions of the existing characterization.
For domain $\P = \Delta_{\Y}$, Theorem \ref{thm:proper-char} formalizes a slightly-informal remark of \citet{gneiting2007strictly}, Section 3.1: given \emph{any} convex function on $\Delta_{\Y}$, \emph{all} of its subtangent rules are proper.
Furthermore, if and only if the function is strictly convex, all of them are strictly proper.
I do not know if proving this result formally would be rigorously possible without much of the development of extended subgradients in this paper.\footnote{This circles us back to first paper to state a characterization, \citet{mccarthy1956measures}, which asserts ``any convex function of a set of probabilities may serve'', but, omitting proofs, merely remarks, ``The derivative has to be taken in a suitable generalized sense.''}

\vskip1em
For general $\P$, the analogous result is \emph{not} true.
Therefore, the question becomes: from which (strictly) convex functions $g$ can we construct a proper subtangent scoring rule?
For such a $g$, does the answer depend on the choice of extended subgradients in (\ref{eqn:subtan})?
Theorem \ref{thm:proper-subset-char} gives a complete answer (see also Figure \ref{fig:scoring-gs}), namely, one must choose certain \emph{$p$-interior-finite} (Definition \ref{def:inter-finite}) extended subgradients of a $g$ that is \emph{interior-locally-Lipschitz} (Definition \ref{def:inter-local}).

\subsection{Tool: extended expected scores}
Before proving the scoring rule characterization, we encounter a problem with the expected score function $S(p;q) := \sum_{y \in \Y} q(y) S(p,y)$ of a regular scoring rule.
Usual characterization proofs proceed by observing that $S(p;q)$ is an affine function of $q$, and that a pointwise supremum over these functions yields the convex $g(q)$ from Definition \ref{def:subtangent-rule} (subtangent rule).

However, $S(p;q)$ is not technically an affine nor affine extended function, as it is only defined on $q \in \Delta_{\Y}$.
Attempting to naively extend it to $q \in \reals^{\Y}$ may lead to illegal sums.\footnote{For example, for the log score, $S(\delta_y,y') = 0$ if $y=y'$, otherwise $-\infty$.
A simple interpretation (e.g. \cite{gneiting2007strictly}) is $S(\delta_y;q) = \sum_{y' \neq y} q(y)(-\infty)$, which raises no contradiction on the simplex, but leads to illegal sums if $q$ contains both positive and negative components.}

Therefore, this section formalizes a key fact about regular scoring rules (Lemma \ref{lemma:exp-score-exist}): for fixed $p$, $S(p;q)$ always coincides on $\Delta_{\Y}$ with at least one affine extended function; we call any such function an extended expected score of $S$.

\begin{definition}[Extended expected score] \label{def:ext-exp-score}
  Let $\Y$ be finite and $\P$ a subset of $\Delta_{\Y}$.
  Given a regular scoring rule $S: \P \times \Y \to \reals \cup \{-\infty\}$, for each $p \in \P$, an affine extended function $h: \reals^{\Y} \to \extreals$ is an \emph{extended expected score function} of $S$ at $p$ if $h(\delta_y) = S(p,y)$ for all $y \in \Y$.
\end{definition}

The name is justified: $h(q)$ is the expected score for report $p$ under belief $q$.
In fact this holds for any $q \in \Delta_{\Y}$, not just $q \in \P$.
Recall that $S(p;q) := \sum_{y \in \Y} q(y) S(p,y)$, Definition \ref{def:proper}.
\begin{observation} \label{obs:ext-exp-score}
  If $h$ is an extended expected score of a regular $S$ at $p$, then for all $q \in \Delta_{\Y}$, $h(q) = S(p;q)$.
\end{observation}
To prove it, write $h(q) = \beta + f(q-p)$ and recall by Definition \ref{def:ext-exp-score} that $S(p,y) = h(\delta_y)$.
So $S(p;q) = \beta + \sum_y q(y) f(\delta_y - p)$.
By regularity of $S$, the sum never contains $+\infty$, so it is a legal sum and, using the scaling axiom of $f$, it equals $\beta + f(q - p)  = h(q)$.

\begin{lemma}[Existence of expected scores] \label{lemma:exp-score-exist}
  Let $\Y$ be finite, $\P \subseteq \Delta_{\Y}$, and $S: \P \times \Y \to \reals \cup \{-\infty\}$ a regular scoring rule.
  Then $S$ has at least one extended expected score function $S_p$ at every $p \in \P$; in fact, it has one that is linear extended.
\end{lemma}
\begin{myproof}
  Given $p$, let $Y_1 = \{y : S(p,y) = -\infty\}$ and $Y_2 = \Y \setminus Y_1$.
  We define $S_p: \reals^{\Y} \to \extreals$ via Algorithm \ref{alg:linext} using the following parameters.
  Set $t = |Y_1|$.
  Let $\hat{f}$ be defined on the subspace spanned by $\{\delta_y : y \in Y_2\}$ via $\hat{f}(x) = \sum_{y \in Y_2} x(y) S(p,y)$.
  By definition, $\hat{f}$ is a finite linear function.
  If $t \geq 1$, let $v_1,\dots,v_t = \{-\delta_y \st y \in Y_1\}$, in any order.
  By Proposition \ref{prop:alg-correct}, $S_p$ is linear extended because it is implemented by Algorithm \ref{alg:linext}.

  To show it is an extended expected score function, we calculate $S_p(\delta_y)$ for any given $y \in \Y$.
  If $y \in Y_2$, then $v_j \cdot \delta_y = 0$ for all $j=1,\dots,t$ and $S_p(\delta_y) = \hat{f}(\delta_y) = S(p,y)$.
  Otherwise, $v_j \cdot \delta_y = 0$ for all $j=1,\dots,t$ except for some $j^*$ where $v_{j^*} = -\delta_y$.
  There, $v_{j^*} \cdot \delta_y = -1$, so $S_p(\delta_y) = -\infty = S(p,y)$.
\end{myproof}

\subsection{Scoring rules on the simplex} \label{sec:scoring-simplex}

This section considers $\P = \Delta_{\Y}$ and uses the machinery of linear extended functions to prove a complete characterization of proper scoring rules, including their construction from any convex function and set of extended subgradients.
This arguably improves little on the state of knowledge as in \citet{gneiting2007strictly} and folklore of the field, but formalizes some important previously-informal or unstated facts.

\begin{theorem}[Construction of proper scoring rules] \label{thm:proper-char}
  Let $\Y$ be a finite set.
  \begin{enumerate}
    \item Let $g: \reals^{\Y} \to \reals \cup \{\infty\}$ be convex with $\effdom{g} \supseteq \Delta_{\Y}$.
          Then: (a) all of its subtangent rules are regular and proper, and there exists at least one; and (b) if $g$ is furthermore strictly convex on $\Delta_{\Y}$, then all of its subtangent rules are strictly proper.
    \item Let $S: \Delta_{\Y} \times \Y \to \reals \cup \{-\infty\}$ be a regular, proper scoring rule.
          Then: (a) it is a subtangent rule of some convex $g$ with $\effdom{g} \supseteq \Delta_{\Y}$, uniquely defined on $\Delta_{\Y}$; and (b) if $S$ is strictly proper, then any such $g$ is strictly convex on $\Delta_{\Y}$; and (c) $S$ is a subgradient rule of some (possibly different) unique positively homogeneous convex $g$.
  \end{enumerate}
\end{theorem}
\begin{myproof}
  (1)
  Given $g$, an arbitrary subtangent rule is defined by $S(p,y) = g(p) + f_p(\delta_y - p)$ for some choice of extended subgradients $f_p: \reals^{\Y} \to \extreals$ at each $p \in \Delta_{\Y}$.
  Proposition \ref{prop:subgrad-exist} asserts that at least one choice of $f_p$ exists for all $p \in \Delta_{\Y}$, so at least one such function $S$ exists.
  And any such $S$ is a regular scoring rule (i.e. never assigns score $+\infty$) because $\effdom{g} \supseteq \Delta_{\Y}$ and by definition of extended subgradient, $S(p,y) \leq g(\delta_y)$.

  To show $S$ is proper: Immediately, the function $S_p(q) = g(p) + f_p(q-p)$ is an extended expected score of $S$ (Definition \ref{def:ext-exp-score}).
  It is an affine extended function supporting $g$ at $p$.
  So $S_q(q) \geq S_p(q)$ for all $p \neq q$, so $S$ is proper by definition, completing \emph{(a)}.
  Furthermore, if $g$ is strictly convex on $\Delta_{\Y}$, then (for any choices of subgradients) each resulting $S_p$ must support $g$ at only one point in $\Delta_{\Y}$ by Lemma \ref{lemma:strict-subgrad}, so $S_q(q) > S_p(q)$ for all $p \neq q$ and $S$ is strictly proper.
  This proves \emph{(b)}.

  $(2)$
  Given $S$, by Lemma \ref{lemma:exp-score-exist}, for each $p$ there exists a linear extended expected score $S_p$.
  Define $g(q) = \sup_{p \in \Delta_{\Y}} S_p(q)$, a convex function by Proposition \ref{prop:sup-char}.
  By definition of properness, $S_q(q) = \max_p S_p(q) = g(q)$.
  So $S_q$ is an affine extended function supporting $g$ at $q \in \effdom{g}$.
  So (e.g. Observation \ref{obs:affine}) it can be written $S_q(x) = g(q) + f_q(x-q)$ for some extended subgradient $f_q$ of $g$ at $q$.
  In particular, $S(p,y) = S(p;\delta_y) = S_p(\delta_y) = g(p) + f_p(\delta_y - p)$, so $S$ is a subtangent rule of $g$, whose effective domain contains $\Delta_{\Y}$.
  Furthermore, by the above, any such $g$ satisfies $g(q) = S_q(q) = S(q;q)$, which is uniquely defined for any $q \in \Delta_{\Y}$.
  This proves \emph{(a)}.

  Furthermore, suppose $S$ is strictly proper and a subtangent rule of some $g$.
  Then as above, it has affine extended expected score functions $S_q$ at each $q$ and, by strict properness, $S_q(p) < S_p(p) = g(p)$ for all $p \in \Delta_{\Y}$, $p \neq q$.
  So at each $p \in \Delta_{\Y}$ there is a supporting affine extended function $S_p$ that supports $g$ nowhere else on $\Delta_{\Y}$.
  By Lemma \ref{lemma:strict-subgrad}, this is a necessary and sufficient condition for strict convexity of $g$ on $\Delta_{\Y}$.
  This proves \emph{(b)}.
  
  Now looking closer, recall Lemma \ref{lemma:exp-score-exist} guarantees existence of a \emph{linear} extended $S_p$.
  So $S_p(\vec{0}) = 0 = g(p) + f_p(-p)$, implying $f_p(p) = g(p) = S_p(p)$.
  Then for all $x$, $S_p(x) = S_p(p) + f_p(x-p)$, which rearranges and uses the axioms of linearity (along with finiteness of $S_p(p)$) to get $S_p(x-p) = f_p(x-p)$, or $S_p = f_p$.
  Since we have $S(p,y) = S_p(\delta_y) = f_p(\delta_y)$, we get that $S$ is a subgradient rule of $g$.
  Furthermore, since each $S_p$ is linear extended, $\alpha g(q) = \alpha \sup_{p \in \Delta_{\Y}} S_p(q) = \sup_{p \in \Delta_{\Y}} \alpha S_p(q) = g(\alpha q)$, for $\alpha \geq 0$.
  So $g$ is positively homogeneous.
  Furthermore, $g(q)$ is uniquely defined on $\Delta_{\Y}$ by $g(q) = S(q;q)$, which then uniquely defines it at any point in the positive orthant.
\end{myproof}

\vskip1em
One nontrivial example is the rule that assigns score $k$ to any prediction with support size $|\Y|-k$, assuming the outcome $y$ is in its support; $-\infty$ if it is not.
This is associated with the convex function $g$ taking the value $k$ on the relative interior of each $|\Y|-k-1$ dimensional face of the simplex; in particular it is zero in the interior of the simplex, $1$ in the relative interior of any facet, \dots, and $|\Y|-1$ at each corner $\delta_{y}$.

A generalization can be obtained from any set function $G: 2^{\Y} \to \reals$ that is \emph{monotone}, meaning $X \subseteq Y \implies G(X) \leq G(Y)$.
We can set $g(p) := G(\Supp{p})$.
In other words, the score for prediction $p$ is $G(\Supp{p})$ if the observation $y$ is in $\Supp{p}$; the score is $-\infty$ otherwise.
We utilized this construction in \citet{chen2016informational}.

Of course neither of these is \emph{strictly} proper.
A similar, strictly proper approach is to construct $g$ from any strictly convex and bounded function $g_0$ on the interior of the simplex; then on the interior of each facet, let $g$ coincide with any bounded strictly convex function whose lower bound exceeds the upper bound of $g_0$; and so on.

%
%
%

\subsection{General belief spaces}

We finally consider general belief and report spaces $\P \subseteq \Delta_{\Y}$, making no restriction except that $\P$ is nonempty.
Here proper scoring rules have been characterized by \citet{gneiting2007strictly} (when $\P$ is convex) and \citet{frongillo2021general} (in general), including for infinite-dimensional outcome spaces.
The statement is that a regular scoring rule $S: \P \times \Y \to \reals \cup \{-\infty\}$ is (strictly) proper if and only if it is a subtangent rule of some (strictly) convex $g$ with\footnote{\citet{frongillo2021general} shows that $S$ is proper on a nonconvex $\P$ only if it coincides with some proper $S'$ on $\text{convhull}(\P)$. This fact will fall out of our result and proof, as $\effdom{g}$ is always convex.} $\effdom{g} \supseteq \P$.

\begin{figure}[ht]
  \caption{Let $\Y = \{0,1\}$; the horizontal axes are $\Delta_{\Y}$ parameterized by $p(1)$. Each subfigure gives the domain of $g$ and whether it is interior-locally-Lipschitz, i.e. produces a scoring rule. \textbf{(\ref{subfig:log})} plots the log scoring rule's associated $g(p) = \sum_y p(y) \log p(y)$; others shift/squeeze/truncate the domain. Vertical supports on the interior of the simplex violate interior-local-Lipschitzness and lead to illegal scoring rules, which would assign e.g. $S(p,0) = +\infty$ and $S(p,1) = -\infty$.}
  \label{fig:scoring-gs}
  \begin{subfigure}{0.24\linewidth}
    \resizebox{\linewidth}{!}{\includegraphics{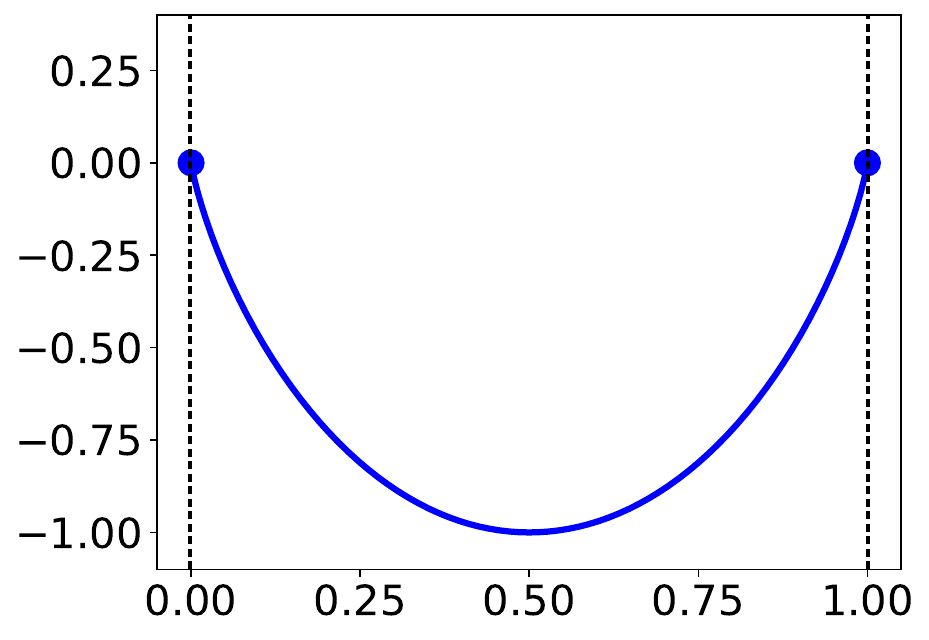}}
    \caption{\footnotesize $[0,1]$; yes.}
    \label{subfig:log}
  \end{subfigure}
  \hfill
  \begin{subfigure}{0.24\linewidth}
    \resizebox{\linewidth}{!}{\includegraphics{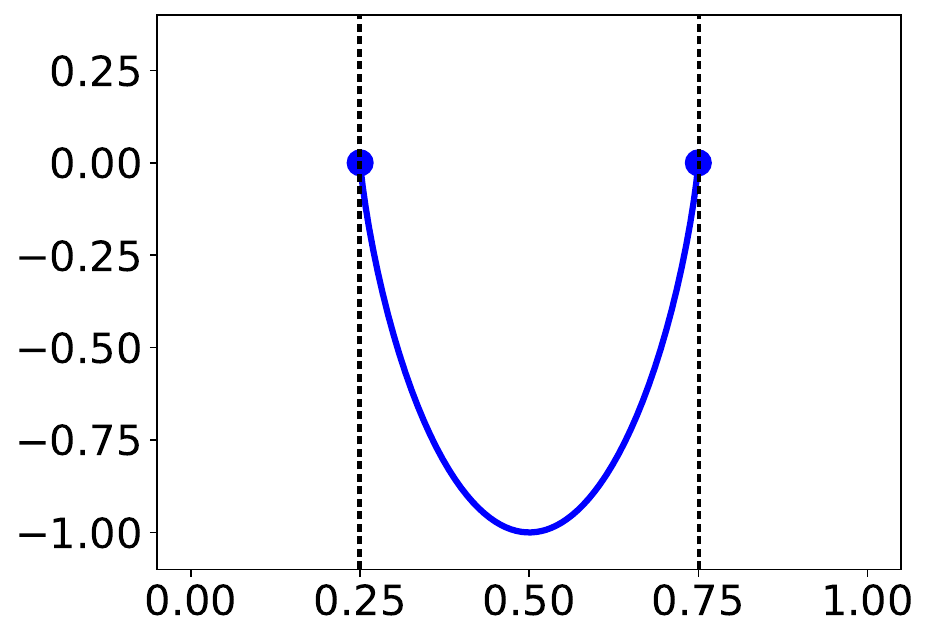}}
    \caption{\footnotesize $[0.25, 0.75]$; no.}
    \label{subfig:shrink-notok}
  \end{subfigure}
  \hfill
  \begin{subfigure}{0.24\linewidth}
    \resizebox{\linewidth}{!}{\includegraphics{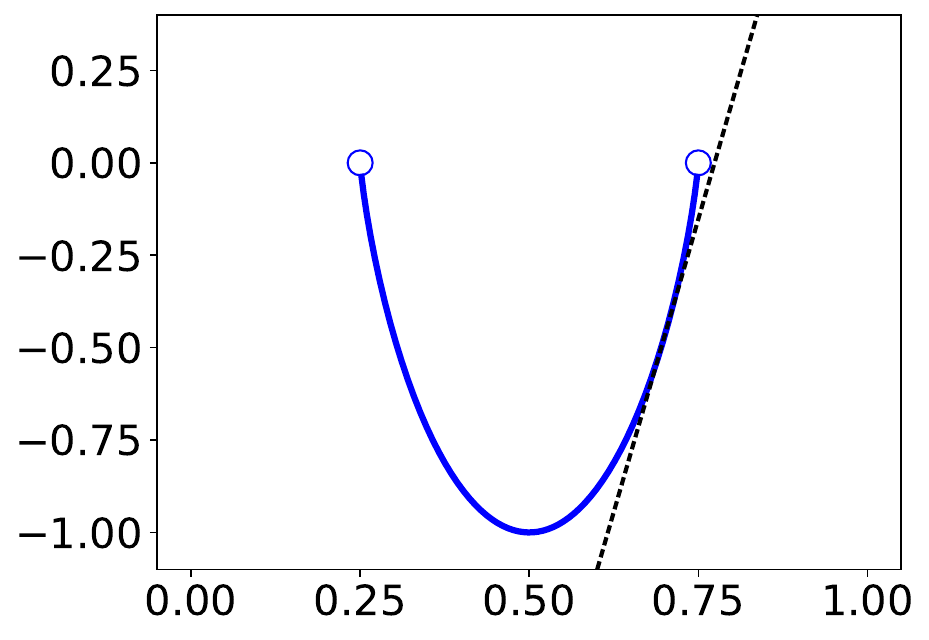}}
    \caption{\footnotesize $(0.25, 0.75)$; yes.}
  \end{subfigure}
  \hfill
  \begin{subfigure}{0.24\linewidth}
    \resizebox{\linewidth}{!}{\includegraphics{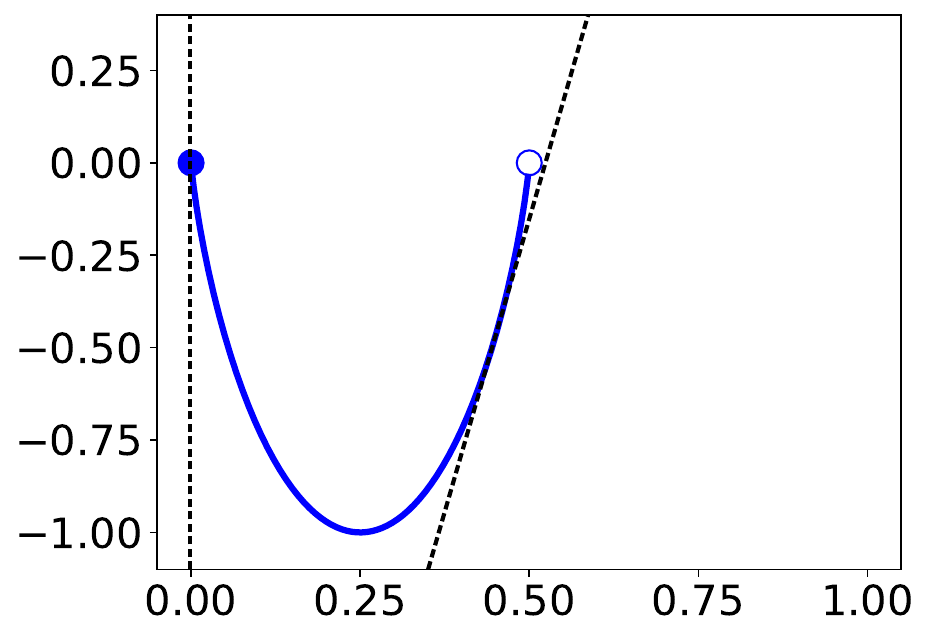}}
    \caption{\footnotesize $[0.0, 0.5)$; yes.}
  \end{subfigure}
\end{figure}

Again, this characterization leaves open the question of exactly which convex functions produce proper scoring rules on $\P$.
Unlike in the case $\P = \Delta_{\Y}$, not all of them do, as Figure \ref{subfig:shrink-notok} illustrates: vertical tangent planes \emph{inside} the simplex lead to non-regular, non-proper scores.

Therefore, the key question is: \textbf{for which convex $g$ are their subtangent rules possibly, or necessarily, regular?}
We next make some definitions to answer this question.

\vskip1em
The following definition will characterize extended subgradients that lead to regular scoring rules.
Recall from Proposition \ref{prop:parsim} that a linear extended function has a unique parameterization in Algorithm \ref{alg:linext}.
Note that because we may have $t=0$ and/or $\Supp{p} = \Y$, the conditions in Definition \ref{def:inter-finite} can be vacuously true.

\begin{definition}[Interior-finite] \label{def:inter-finite}
  Say a linear extended function $f: \reals^{\Y} \to \extreals$ is \emph{$p$-interior-finite} for $p \in \Delta_{\Y}$ if its parameterization $(t,v_1,\dots,v_t,\hat{f})$ in Algorithm \ref{alg:linext} has: (1) for all pairs $y,y' \in \Supp{p}$, that $v_j(y) = v_j(y')$ for all $j=1,\dots,t$; and (2) for all $y \not\in \Supp{p}$ and $y' \in \Supp{p}$, that $v_j(y) = v_j(y')$ for some sequence $j=1,\dots,k$, with either $k=t$ or else $v_{k+1}(y) < v_{k+1}(y')$.
\end{definition}

If $f$ is $p$-interior-finite, then property (1) gives that $f(q-p)$ is finite for $q$ in the face of the simplex generated by the support of $p$.
Property (2) gives that $f(q-p)$ is either finite or $-\infty$ at all $q \in \Delta_{\Y}$ that are not in that face.
These are the two key properties that we will need for regular scoring rules.

\begin{lemma}[Interior-finite implies quasi-integrable] \label{lemma:interior-quasi}
  Let $p \in \Delta_{\Y}$.
  A linear extended function $f: \reals^{\Y} \to \extreals$ is $p$-interior-finite if and only if $f(\delta_y - p) \in \reals \cup \{-\infty\}$ for all $y \in \Y$. 
  Furthermore, in this case $f(\delta_y - p) \in \reals$ for any $y \in \Supp{p}$.
\end{lemma}
\begin{myproof}
  Let $f$ be a linear extended function and $v_1,\dots,v_t$ its unit vectors in the parameterization of Algorithm \ref{alg:linext}.

  $(\implies)$
  Suppose $f$ is $p$-interior-finite.
  If it is finite everywhere, then the conclusion is immediate.
  Otherwise, let $y \in \Supp{p}$.
  Using that $v_1(y) = v_1(y')$ for any $y' \in \Supp{p}$, we have
  \begin{align*}
    v_1 \cdot (\delta_y - p) &= v_1(y) - v_1 \cdot p  \\
                             &= v_1(y) - \sum_{y' \in \Supp{p}} v_1(y') p(y')  \\
                             &= v_1(y) - v_1(y) \sum_{y' \in \Supp{p}} p(y')  \\
                             &= 0 .
  \end{align*}
  This applies to each successive $v_j$, so $f(\delta_y - p) \in \reals$.

  Now let $y \not\in \Supp{p}$.
  Having just shown $v_1 \cdot (\delta_{y'} - p) = 0$ for any $y' \in \Supp{p}$, we have $v_1 \cdot (\delta_y - p) = v_1 \cdot (\delta_y - \delta_{y'}) = v_1(y) - v_1(y') \leq 0$ by assumption of $p$-interior-finite.
  So either Algorithm \ref{alg:linext} returns $-\infty$ in the first iteration, or we continue to the second iteration.
  The same argument applies at each iteration, so $f(\delta_y - p) \in \reals \cup \{-\infty\}$.

  $(\impliedby)$
  Suppose $f(\delta_y - p) \in \reals \cup \{-\infty\}$ for all $y \in \Y$.
  Again if $f$ is finite everywhere, it is $p$-interior-finite, QED.
  So suppose its rank is at least $1$.
  We must have $v_1 \cdot (\delta_y - p) \leq 0$ for all $y \in \Y$.
  This rearranges to $v_1(y) \leq v_1 \cdot p$ for all $y$, or $\max_{y \in \Y} v_1(y) \leq v_1 \cdot p$.
  But $p$ is a probability distribution, so $v_1 \cdot p \leq \max_{y \in \Y} v_1(y)$.
  So $v_1 \cdot p = \max_{y \in \Y} v_1(y)$.
  This implies that $v_1(y) = v_1(y') = v_1 \cdot p$ for all $y,y' \in \Supp{p}$.
  So $v_1 \cdot (\delta_y - p) = 0$ if $y \in \Supp{p}$.
  Meanwhile, if $y \not\in \Supp{p}$, then $v_1(y) \leq v_1 \cdot p = v_1(y')$ for any $y' \in \Supp{p}$.

  This argument repeats to prove that if $y \in \Supp{p}$, then for all $j=1,\dots,t$, $v_j(y) = \max_{y' \in \Y} v_j(y')$ and also $v_j \cdot (\delta_y - p) = 0$.
  So it also gives $f(\delta_y - p) \in \reals$.
  Meanwhile, if $y \not\in \Supp{p}$, then we have a sequence $v_1(y) = \max_{y' \in \Y} v_1(y')$, \dots, $v_k(y) = \max_{y' \in \Y} v_k(y')$, where in each case $v_j \cdot (\delta_y - p) = 0$.
  Then finally, either $k=t$ and $f(\delta_y - p) \in \reals$, or we have $v_{k+1}(y) < \max_{y' \in \Y} v_{k+1}(y')$ and $f(\delta_y - p) = -\infty$.
\end{myproof}

We can now make the key definitions for the scoring rule characterization.
\begin{definition}[Interior-locally-Lipschitz; ILL subtangent rule] \label{def:inter-local}
  For $\P \subseteq \Delta_{\Y} \cap \effdom{g}$, say $g: \reals^{d} \to \reals \cup \{\infty\}$ is \emph{$\P$-interior-locally-Lipschitz} if $g$ has at every $p \in \P$ a $p$-interior-finite extended subgradient.
  Say that a scoring rule $S: \P \times \Y \to \extreals$ is an \emph{ILL subtangent rule of $g$} if it is a subtangent rule of $g$ (Definition \ref{def:subtangent-rule}) where each extended subgradient $f_p$ is $p$-interior finite.
\end{definition}

An interior-locally-Lipschitz $g$ simply has two conditions on its vertical supporting hyperplanes.
First, they cannot cut faces of the simplex whose relative interior intersects $\effdom{g}$; they can only contain them.
This is enforced by $v_j(y) = v_j(y')$ for $y,y' \in \Supp{p}$.
In other words, the extended subgradients must be \emph{finite} (in particular bounded, an analogy of causing $g$ to be Lipschitz) relative to the affine hull of these faces.
Second, they must be oriented correctly so as not to cut that face from the rest of the simplex.
This is enforced by $v_j(y) \leq v_j(y')$ for $y \not\in \Supp{p}$, $y' \in \Supp{p}$.

Lemma \ref{lemma:interior-quasi} gives the following corollary.
\begin{corollary} \label{cor:lip-reg}
  A subtangent scoring rule $S: \P \times \Y \to \extreals$ is regular if and only if it is an ILL subtangent rule.
  A convex $g: \reals^{\Y} \to \reals \cup \{\infty\}$ with $\effdom{g} \supseteq \P$ has a regular subtangent scoring rule if and only if $g$ is $\P$-interior-locally-Lipschitz.
\end{corollary}
\begin{myproof}
  Let $S$ be a subtangent scoring rule of some $g$ with $\effdom{g} \supseteq \P$, with extended subgradients $\{f_p : p \in \P\}$
  Recall that $S$ is regular if $S(p,y) \in \reals \cup \{-\infty\}$ for all $p \in \P$ and $y \in \Y$.
  Because $g(p) \in \reals$ for all $p \in \P$, regularity is equivalent to $f_p(\delta_y - p) \in \reals \cup \{-\infty\}$ for all $p \in \P$ and all $y \in \Y$.
  This is equivalent to $p$-interior-local-Lipschitzness of each $f_p$ by Lemma \ref{lemma:interior-quasi}, proving the first claim.
  By definition, $g$ is $\P$-interior-locally-Lipschitz if and only if it has an ILL subtangent rule, proving the second claim.
\end{myproof}

\begin{theorem}[Construction of scoring rules on $\P$] \label{thm:proper-subset-char}
  Let $\Y$ be a finite set and $\P$ an arbitrary nonempty subset of $\Delta_{\Y}$.
  \begin{enumerate}
    \item Let $g: \reals^{\Y} \to \reals \cup \{\infty\}$ be convex and $\P$-interior-locally-Lipschitz with $\effdom{g} \supseteq \P$.
          Then: (a) a subtangent scoring rule of $g$ is regular and proper if and only if it is an ILL subtangent rule of $g$, and there exists at least one; and (b) if $g$ is strictly convex on $\P$, then all of its ILL subtangent rules are regular and strictly proper.
    \item Let $S: \P \times \Y \to \reals \cup \{-\infty\}$ be a regular proper scoring rule.
          Then: (a) it is an ILL subtangent rule of some convex $\P$-interior-locally-Lipschitz $g: \reals^{\Y} \to \reals \cup \{\infty\}$ with $\effdom{g} \supseteq \P$, uniquely defined on $\P$; and (b) if $S$ is strictly proper and $\P$ a convex set, then any such $g$ is strictly convex on $\P$.  
  \end{enumerate}
\end{theorem}
Observe that in (1), $g$ may have some subtangent rules that are regular and some that are not, depending on the choices of extended subgradients at each $p$.
An example is if $\P = \effdom{g}$ is a single point, say the uniform distribution; vertical tangent planes will not work, but finite subgradients will.
(This issue did not arise in the previous section where $\P = \Delta_{\Y}$; there, all of a convex $g$'s subtangent rules were regular.)
\begin{myproof}
  (1)
  Let such a $g$ be given.
  The definition of $\P$-interior-locally-Lipschitz immediately implies that it has at least one ILL subtangent rule.
  Let $S: \P \times \Y \to \extreals$ be any such rule.
  Then $S$ is of the form $S(p,y) = g(p) + f_p(\delta_y - p)$.
  Applying Lemma \ref{lemma:interior-quasi}, $f_p(\delta_y - p) \in \reals \cup \{-\infty\}$, so $S(p,y) \neq \infty$, so it is a regular scoring rule.
  Immediately we have an affine extended expected score function $S_p$ at each $p$ satisfying $S_p(q) = g(p) + f_p(q-p)$ for all $p \in \P$ and all $q \in \Delta_{\Y}$.
  We have $S_p(p) \geq S_q(p)$ for all $q$ because $S_p$ supports $g$ at $p$, so $S$ is proper.

  Conversely, if $S$ is a subtangent rule but not an ILL one, then there is some $p \in \P$ where $S(p,y) = g(p) + f_p(\delta_y - p)$ where $f_p$ is not $p$-interior-finite.
  By Lemma \ref{lemma:interior-quasi}, there exists $y$ such that $S(p,y) = \infty$, so $S$ is not regular.
  (Recall that by our definition, only regular scoring rules can be proper.)
  This proves \emph{(a)}.

  Now suppose $g$ is strictly convex and $S$ an ILL subtanget rule.
  Then each $S_q$ supports $g$ at just one point in its effective domain (Lemma \ref{lemma:strict-subgrad}).
  So $S_p(p) > S_q(p)$ for $p,q \in \P$ with $q \neq p$, so $S$ is strictly proper, proving \emph{(b)}.

  (2)
  Let a regular proper scoring rule $S$ be given.
  It has by Lemma \ref{lemma:exp-score-exist} an affine extended expected score function $S_p$ at each $p \in \P$.
  Define $g: \reals^{\Y} \to \reals \cup \{\infty\}$ by $g(q) = \sup_{p \in \P} S_p(q)$.
  Then $g$ is convex by Proposition \ref{prop:sup-char}.
  By definition of proper, for all $p \in \P$, we have $S_p(p) = \max_{q \in \P} S_q(p) = g(p)$.
  So $S_p$ is an affine extended function supporting $g$ at $p$, so it can be written $S_p(q) = g(p) + f_p(q-p)$ for some subgradient $f_p$ of $g$ at $p$.
  By definition of extended expected score, $S(p,y) = S_p(\delta_y) = g(p) + f_p(\delta_y - p)$, so it is a subtangent rule.

  Next, since $S$ is a scoring rule, in particular $S_p(\delta_y) \in \reals \cup \{-\infty\}$ for all $y$.
  So by Lemma \ref{lemma:interior-quasi}, each $f_p$ is $p$-interior-finite, so $g$ is interior-locally-Lipschitz and $S$ is its ILL subtangent rule.
  Finally, observe that if $S$ is a subtangent rule of any $g$, then $g(p) = S_p(p)$ for any $p \in \P$, so $g$ is uniquely defined on $\P$.
  This completes the proof of \emph{(a)}.

  Now suppose $S$ is strictly proper and $\P$ a convex set.
  Then for any $q \in \P$, we have that $S_q$ supports $g$ at $q$ alone out of all points in $\P$, by definition of strict properness.
  By Lemma \ref{lemma:strict-subgrad}, $g$ is strictly convex on $\P$, proving \emph{(b)}.
\end{myproof}

There remain open directions investigating interior-locally-Lipschitz functions.
Two important cases can be immediately resolved:
\begin{itemize}
  \item If $\P \subseteq \inter{\Delta_{\Y}}$, then $g$ is $\P$-interior-locally-Lipschitz if $g$ is subdifferentiable on $\P$, i.e. has finite subgradients everywhere on $\P$.
        This follows because every $p \in \P$ has full support.
        And in fact, it is necessary and sufficient that $g$ be subdifferentiable when restricted to the affine hull of the simplex; its subgradients can still be infinite perpendicular to the simplex, e.g. of the form $f(x) = \infty \cdot (x \cdot \vec{1}) + \hat{f}(x)$.
  \item If $\P$ is an open set relative to the affine hull of $\Delta_{\Y}$, then every convex function $g$ with $\effdom{g} \supseteq \P$ is $\P$-interior-locally-Lipschitz.
        This follows from the previous case along with the fact that convex functions are subdifferentiable on open sets in their domain.
\end{itemize}
One question for future work is whether one can optimize efficiently over the space of interior-locally-Lipschitz $g$, for a given set $\P \subseteq \Delta_{\Y}$.
This could be useful in constructing proper scoring rules that optimize some objective.

Another direction is to extend these constructions to general affine type spaces as studied by \citet{frongillo2021general}, where the structure of the probability simplex is no longer present.

\bibliographystyle{plainnat}
\bibliography{citations}

\end{document}